\documentclass[11pt]{amsart}

\usepackage{amsmath, amssymb, amscd}
\usepackage{eucal}
\usepackage{mathrsfs}
\usepackage[frame,cmtip,arrow,matrix,line,graph,curve]{xy}
\usepackage{graphicx}
\usepackage[colorlinks=true,citecolor=red,linkcolor=blue]{hyper ref}

\usepackage{pgf,tikz}
\usetikzlibrary{arrows}
\usetikzlibrary{decorations.pathmorphing}
\usetikzlibrary{decorations.pathreplacing,decorations.markings}
\tikzset{arrow data/.style 2 args={%
      decoration={%
         markings,
         mark=at position #1 with \arrow{#2}},
         postaction=decorate}
      }%

\definecolor{zzttqq}{rgb}{0.6,0.2,0.}
\definecolor{qqqqff}{rgb}{0.,0.,1.}

\newcommand{\cC}{\mathcal C}

\newcommand{\oD}{\overline{D}}
\newcommand{\oGamma}{\overline{\Gamma}}
\newcommand{\oPhi}{\overline{\Phi}}
\newcommand{\oO}{\overline{O}}
\newcommand{\oS}{\overline{S}}

\newcommand{\tC}{{\widetilde C}}

\newcommand{\tgamma}{{\widetilde{\gamma}}}
\newcommand{\tGamma}{{\widetilde{\Gamma}}}

\newcommand{\tPhi}{\widetilde{\Phi}}

\newcommand{\tS}{\widetilde S}
\newcommand{\tv}{\widetilde{v}}
\newcommand{\tx}{\widetilde{x}}

\newcommand{\rr}{\mathbb R}
\newcommand{\zz}{\mathbb Z}

\newcommand{\ep}{\epsilon}

\newcommand{\cl}{\mathrm{cl} }
\newcommand{\Bd}{\mathrm{Bd} }

\newcommand{\Int}{\mathrm{Int}}
\newcommand{\inv}{^{-1}}
\newcommand{\op}{\mathrm{op} }
\newcommand{\psr}{\mathrm{psr} }
\newcommand{\sr}{\mathrm{sr} }
\newcommand{\sing}{\mathrm{sing}}
\newcommand{\Sing}{\mathrm{Sing}}
\newcommand{\sm}{\mathrm{sm}}

\newtheorem{prop}{Proposition}[section]
\newtheorem{thm}[prop]{Theorem}
\newtheorem{lem}[prop]{Lemma}
\newtheorem{cor}[prop]{Corollary}

\theoremstyle{remark}
  \newtheorem{rk}[prop]{Remark}
  
\theoremstyle{definition}
 \newtheorem{example}[prop]{Example}
 \newtheorem{defn}[prop]{Definition}

\numberwithin{equation}{section}
\numberwithin{figure}{section}

\begin{document}
\title[Polycycle omega-limit sets of flows on the compact Riemann surfaces and Eulerian path]
{Polycycle omega-limit sets of flows on the compact Riemann surfaces and Eulerian path}
\author{Jaeyoo Choy$^{\dagger}$ and Hahng-Yun Chu$^{\ast}$}
\address{Dept.\ Mathematics, Kyungpook National University,
Sankyuk-dong, Buk-gu,
Daegu 702-701, Republic of Korea}\email{\rm choy@knu.ac.kr\ (J. Choy)}
\address{Dept.\ Mathematics, Chungnam National
University, 79, Daehangno, Yuseong-gu,
Daejeon 305-764, Republic of Korea}
\email{\rm hychu@cnu.ac.kr\ (H.-Y.Chu)}

\thanks{$\dagger$  The first author is partially supported by Kyungpook National University Fund 2017.}

\thanks{$\ast$ Corresponding author. His research
is  performed as a subproject of project Research for Applications of Mathematical Principles (No C21501) and supported by the National Institute of Mathematics Sciences(NIMS)}


\begin{abstract}
Let $(S,\Phi)$ be a pair of a closed oriented surface and $\Phi$ be a real analytic flow with finitely many singularities.
Let $x$ be a point of $S$ with the polycycle $\omega$-limit set $\omega(x)$.
In this paper we give topological classification of $\omega(x)$.
Our main theorem says that $\omega(x)$ is diffeomorphic to the boundary of a cactus in the $2$-sphere $S^{2}$.
Moreover $S$ is a connected sum of the above $S^{2}$ and a closed oriented surface along finitely many embedded circles which are disjoint from $\omega(x)$.
This gives a natural generalization to the higher genus of the main result of \cite{JL} for the genus $0$ case.

Our result is further applicable to a larger class of surface flows, a compact oriented surface with corner and a $\cC^{1}$-flow with finitely many singularities locally diffeomorphic to an analytic flow.
\end{abstract}

\maketitle

\thispagestyle{empty} \markboth{Jaeyoo Choy and Hahng-Yun
Chu}{Polycycle omega-limit sets of flows on the compact Riemann surfaces and Eulerian path}

\section{Introduction}\label{sec: intro}


\subsection{Main theorem}\label{subsec: main theorem}

Let $M$ be a $\cC^{1}$-manifold with a $\cC^{1}$-flow $\Phi\colon \rr\times M\to M$.
For $x\in M$, we denote the $\omega$-limit set of $x$ by $\omega(x)$.

A main purpose of the paper is to give complete topological
classification of $\omega$-limits for (real)
analytic flows on closed oriented surface (= Riemann surface).
The Riemann sphere case was done in Jim\'enez L\'opez-Llibre's work \cite{JL} (cf.\ \cite{BL}).
In this case, $\omega(x)$ is either a single point or the boundary of a cactus.
Here a cactus is a simply connected finite union of disks mutually intersecting at most one point (Fig.\ \ref{Fig: 1-A}; see the definition in Example \ref{example: semi ribbon} (1)).
\begin{figure}[ht]
\begin{center}
\begin{tikzpicture}[line cap=round,line join=round,>=triangle 45,x=0.7cm,y=0.7cm]
\draw plot [smooth, tension=.2] coordinates
{(4.46,0.34) (3.2,1.68) (3.12,2.98) (3.84,3.92) (5.14,3.94) (5.88,3.02) (5.8,1.78) (4.46,0.34)};
\fill[color=zzttqq,fill=zzttqq,fill opacity=0.1] plot [smooth, tension=.2] coordinates
{(4.46,0.34) (3.2,1.68) (3.12,2.98) (3.84,3.92) (5.14,3.94) (5.88,3.02) (5.8,1.78) (4.46,0.34)};
\draw plot [smooth, tension=.2] coordinates
{(4.46,0.34) (2.88,0.4) (1.8,-0.36)   (1.74,-1.68)   (3.1,-1.94)   (4.2,-1.3) (4.46,0.34)};
\fill[color=zzttqq,fill=zzttqq,fill opacity=0.1] plot [smooth, tension=.2] coordinates
{(4.46,0.34)   (2.88,0.4)   (1.8,-0.36)   (1.74,-1.68)   (3.1,-1.94)   (4.2,-1.3) (4.46,0.34)};
\draw plot [smooth, tension=.2] coordinates
{(4.46,0.34)   (5.98,0.58)   (7.42,-0.08)   (7.52,-1.32)   (6.38,-1.84)   (5.18,-1.1) (4.46,0.34)};
\fill[color=zzttqq,fill=zzttqq,fill opacity=0.1] plot [smooth, tension=.2] coordinates
{(4.46,0.34)   (5.98,0.58)   (7.42,-0.08)   (7.52,-1.32)   (6.38,-1.84)   (5.18,-1.1) (4.46,0.34)};
\draw plot [smooth, tension=.2] coordinates
{(7.42,-0.08)   (7.78,1.04)   (8.9,1.48)   (9.86,1.32)   (9.52,0.42)   (8.54,-0.06) (7.42,-0.08)};
\fill[color=zzttqq,fill=zzttqq,fill opacity=0.1] plot [smooth, tension=.2] coordinates
{(7.42,-0.08)   (7.78,1.04)   (8.9,1.48)   (9.86,1.32)   (9.52,0.42)   (8.54,-0.06) (7.42,-0.08)};
\draw plot [smooth, tension=.2] coordinates
{(1.8,-0.36)   (1.6,0.54)   (0.74,0.86)   (0.04,0.44)   (0.24,-0.22)   (1.04,-0.5) (1.8,-0.36)};
\fill[color=zzttqq,fill=zzttqq,fill opacity=0.1] plot [smooth, tension=.2] coordinates
{(1.8,-0.36)   (1.6,0.54)   (0.74,0.86)   (0.04,0.44)   (0.24,-0.22)   (1.04,-0.5) (1.8,-0.36)};
\draw plot [smooth, tension=.2] coordinates
{(9.52,0.42)   (9.34,-0.36)   (9.66,-0.9)   (10.3,-0.92)   (10.7,-0.36)   (10.34,0.12) (9.52,0.42)};
\fill[color=zzttqq,fill=zzttqq,fill opacity=0.1] plot [smooth, tension=.2] coordinates
{(9.52,0.42)   (9.34,-0.36)   (9.66,-0.9)   (10.3,-0.92)   (10.7,-0.36)   (10.34,0.12) (9.52,0.42)};
\draw plot [smooth, tension=.2] coordinates
{(8.9,1.48)   (9.2,2.16)   (8.9,2.8)   (8.18,2.88)   (7.94,2.24)   (8.12,1.66) (8.9,1.48)};
\fill[color=zzttqq,fill=zzttqq,fill opacity=0.1] plot [smooth, tension=.2] coordinates
{(8.9,1.48)   (9.2,2.16)   (8.9,2.8)   (8.18,2.88)   (7.94,2.24)   (8.12,1.66) (8.9,1.48)};
\end{tikzpicture}
\end{center}
\caption{a cactus}
\label{Fig: 1-A}
\end{figure}

A version of our result generalizes a result of Jim\'enez L\'opez-Llibre \cite[Theorem C]{JL} on $S^{2}$ to an arbitrary Riemann surface.

\begin{thm}\label{th: main: analytic version}
Let $(S,\Phi)$ be a local real analytic flow on a Riemann surface with only finitely many singularities.
Let $x\in S$.
Suppose the $\omega$-limit set $\omega(x)$ is a polycycle.
Then we have the following description of $S$ and $\Phi$.
	\begin{itemize}
	\item[(i)] $S$ is a connected sum of $S^{2}$ and another Riemann surface along finitely many embedded circles disjoint from $\omega(x)$.
	\item[(ii)] There is a cactus $D$ in $S^{2}$ such that its boundary $\Bd(D)$ coincides with $\omega(x)$. 
	\item[(iii)]
Moreover there is an open neighborhood $U$ of $D$ in $S^{2}$ such that $\Phi(t,x)\in U\setminus D$ for all sufficiently large $t$ and  $\lim_{t\to\infty}d(\Phi(t,y),D)=0$ for each $y\in U\setminus D$. 
	\end{itemize}
\end{thm}

Here a \textit{polycycle} is defined as follows: Let $O_{1},O_{2},...,O_{n}$ be mutually disjoint free orbits.
We suppose that there are (\textit{not} necessarily distinct) singular points $x_{1},x_{2},...,x_{n}$ such that $\alpha(z)=\{x_{l}\},\ \omega(z)=\{x_{l+1}\}$ for $z\in O_{l}$ and each $l\in\{1,2,...,n\}$, where $x_{n+1}$ is set to be $x_{1}$. 
Here a \textit{free} orbit means a non-periodic non-singular orbit.
Now a polycycle is defined to be the union of the closures $\oO_{1}\cup\oO_{2}\cup...\cup\oO_{n}$.

\begin{figure}[ht]
\begin{center}
\begin{tikzpicture}[line cap=round,line join=round,>=triangle 45,x=0.6cm,y=0.6cm]
\draw [rotate around={0.:(6.8,-3.6)}] (6.8,-3.6) ellipse (7.4726166769077675 and 2.);
\draw plot[smooth, tension=.2pt] coordinates
{(3.2,-2.8) (2.4,-3.6) (3.2,-4.4) (4.4,-4.4) (5.2,-3.6) (4.34,-2.8) (3.2,-2.8)};
\draw plot[smooth, tension=.2pt] coordinates
{(5.2,-3.6) (6.,-2.8) (7.2,-2.8) (8.,-3.6) (7.2,-4.4) (6.,-4.4) (5.2,-3.6)};
\draw [line width=1.5pt](3.8,-3.6) circle (0.6);
\draw [line width=1.5pt](6.6,-3.58) circle (0.6003332407921459);
\draw [line width=1.5pt](9.42,-3.6) circle (0.58);
\draw [line width=1.5pt](9.42,-0.42) circle (0.5803447251418764);
\draw [line width=1.5pt](6.62,-0.42) circle (0.5803447251418764);
\draw [line width=1.5pt] (3.8,-0.42) circle (0.6003332407921459);
\draw plot[smooth, tension=.6pt] coordinates
{(7.2,-0.4) (7.6,0.8) (7.98,1.) (8.4,0.8) (8.839956987758756,-0.401288935088992)};
\draw plot[smooth, tension=.6pt] coordinates
{(10.,-0.4) (10.4,1.2) (11.2,1.6) (12.8,2.) (13.6,2.8) (13.6,3.6) (12.8,4.4) (11.2,4.8) (8.8,4.8) (7.2,4.) (6.4,2.8) (6.039956987758756,-0.40128893508899227)};
\draw [rotate around={0.:(8.4,2.8)}] (8.4,2.8) ellipse (0.8944271909998783 and 0.4);
\draw [rotate around={0.:(11.2,2.8)}] (11.2,2.8) ellipse (0.8944271909998164 and 0.4);
\draw plot[smooth, tension=.6pt] coordinates
{ (3.2,-0.4) (2.8,1.14) (2.,1.6)  (0.4,2.)  (-0.4,2.8)  (-0.4,3.6) (0.4,4.4) (2.,4.8) (3.2,4.8)  (4.4,4.)  (4.8,2.8) (4.8,1.6) (4.4,-0.4)};
\draw [rotate around={0.:(2.4,3.2)}] (2.4,3.2) ellipse (0.8944271909999166 and 0.4);
\end{tikzpicture}
\end{center}
\caption{connected sum of $S^2$ and a (disconnected) Riemann surface of genus $3$ along the bold circles; ellipsoids stand for genii of Riemann surface; two hexagons form an $\omega$-limit set}
\label{Fig: 1-B}
\end{figure}
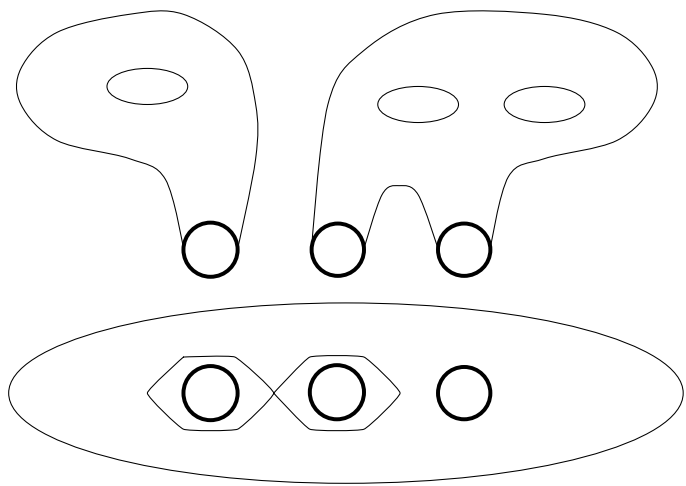

The connected sum and $\omega(x)$ in the statement of the theorem are depicted as in Fig.\ \ref{Fig: 1-B}.

This theorem will be reformulated for wider classes of surfaces, namely surfaces with corner in Theorem \ref{th: main: omega}.
Surface with corner is more natural object for  $\cC^{1}$-flows because cutting a closed  $\cC^{1}$-surface along an invariant simple closed curve gives a compact  $\cC^{1}$-surface with corner.
In Theorem \ref{th: main: omega} we also obtain an explicit example of $\cC^{1}$-flows on a given oriented surface with corner and a prescribed cactus boundary described in (i)--(iii) above.
The procedure involves topological surgery, so does \textit{not} assure construction of real analytic flows, which we leave it as a question. 

Let us give further explication on the theorem.
First we look back the polycycle assumption in the statement.
If $S$ has genus $\ge1$, the $\omega$-limit set $\omega(x)$ is \textit{not} necessarily a finite union of orbits in general.
E.g., on the $2$-torus $T^{2}$ with a flow with an irrational slope, any $\omega$-limit set becomes the whole space $T^{2}$, an infinite union of orbits.
This is a new phenomenon which does not happen in genus $0$.
However the $\omega$-limit set being polycycle is \textit{not} just technical assumption. 
For instance, it is a consequence in theorems of Palis-de Melo \cite{PM} and Bendixson \cite{Be}, which relate it with the Anosov property or trivial recurrence property (see Corollary \ref{cor: Pujals de Melo Bendixson} for the precise assertion).

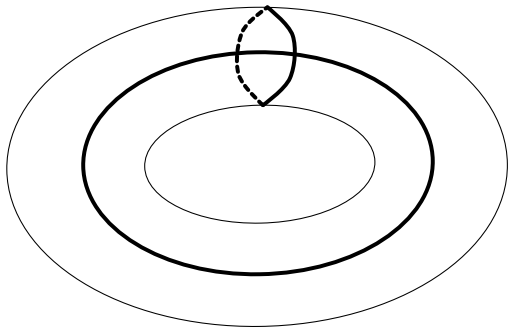
\begin{figure}[ht]
\begin{center}
\begin{tikzpicture}[line cap=round,line join=round,>=triangle 45,x=0.6cm,y=0.6cm]
\draw [rotate around={0.6708801057280257:(5.35,0.27)}] (5.35,0.27) ellipse (5.5459475529580295 and 3.538662778531059 );
\draw [rotate around={1.3078956016637846:(5.41,0.33)}] (5.41,0.33) ellipse (2.5510369357888405  and 1.307359723931758 );
\draw [rotate around={0.9580341086759661:(5.37,0.35)},line width=1.5pt] (5.37,0.35) ellipse (3.8729949955402954 and 2.4611969111552567 );
\draw[line width=1.5pt] plot [smooth, tension=.6, line width=1.5pt] coordinates {(5.579753367429981,3.8073633057197767) (6.14,3.18)  (6.08,2.24) (5.480168006653967,1.638297082190812)};
\draw[dash pattern=on 2pt off 2pt, line width=1.5pt]  plot [smooth, tension=.6, line width=1.5pt] coordinates
{(5.480168006653967,1.638297082190812) (4.94,2.32) (5.02,3.26) (5.579753367429981,3.8073633057197767)};
\end{tikzpicture}
\end{center}
\caption{the union of a latitude and a longitude on $T^2$}
\label{Fig: 1-C}
\end{figure}

Secondly a polycyle $\oO_{1}\cup\oO_{2}\cup...\cup\oO_{n}$ with the singular points $x_{1},x_{2},...,x_{n}$ is realized as an embedded finite oriented graph by the correspondence 
	$$
	O_{l}\leftrightarrow \mbox{an edge}, \quad x_{l}\leftrightarrow\mbox{a vertex} 
	$$
(see \S\ref{subsec: semiribbon} for further explanation).
We denote this graph by $\Gamma$.
However \textit{not} every finite oriented graph is realized as an $\omega$-limit set, even if it is isomorphic to the boundary of a cactus on $S^{2}$ (a \textit{cactus boundary} for short).
For, Theorem \ref{th: main: analytic version} gives a strong regulation on topology of $(S,\omega(x))$.
For instance $\Gamma$ is necessarily a hyperbolic graph. 
Thus the following example cannot be an $\omega$-limit set for any analytic flow with finitely many singularities: the union of a latitude and a longitude in $T^{2}$  (Fig.\ \ref{Fig: 1-C}).

\begin{figure}[ht]
\begin{center}
\begin{tikzpicture}[line cap=round,line join=round,>=triangle 45,x=0.6cm,y=0.6cm]
\draw(0.4,-0.58) circle (3.58);
\draw (7.78,2.9)-- (7.82,-3.84);
\draw (7.82,-3.84)-- (14.56,-3.8);
\draw (14.56,-3.8)-- (14.52,2.94);
\draw (14.52,2.94)-- (7.78,2.9);
\draw [line width=1.2pt] (0.4,-0.58)-- (1.58,0.54);
\draw [line width=1.2pt] (1.58,0.54)-- (2.86,-0.4);
\draw [line width=1.2pt] (2.86,-0.4)-- (1.64,-1.58);
\draw [line width=1.2pt] (1.64,-1.58)-- (0.4,-0.58);
\draw [line width=1.2pt] (0.4,-0.58)-- (-0.7,0.42);
\draw [line width=1.2pt] (-0.7,0.42)-- (-1.86,-0.68);
\draw [line width=1.2pt] (-1.86,-0.68)-- (-0.72,-1.78);
\draw [line width=1.2pt] (-0.72,-1.78)-- (0.4,-0.58);
\draw[line width=1.2pt, color=gray] plot [smooth, tension=.3pt] coordinates {(-0.92,0.72) (-2.18,-0.46) (-2.2,-1.02) (-1.08,-2.06) (-0.48,-2.12)  (0.22,-1.36) (0.7,-1.34) (1.44,-1.96) (1.92,-1.98) (3.18,-0.7) (3.2,-0.12) (1.92,0.84) (1.32,0.8) (0.68,0.14) (0.24,0.1) (-0.38,0.72) (-0.92,0.72)};
\draw (7.78,2.9)-- (7.82,-3.84);
\draw (7.8012462688540385,-0.6799963019053985) -- (7.665625511773548,-0.5757993237874383);
\draw (7.8012462688540385,-0.6799963019053985) -- (7.93562075708049,-0.5741969781179603);
\draw (7.8,-0.47) -- (7.664379242919511,-0.3658030218820397);
\draw (7.8,-0.47) -- (7.934374488226451,-0.3642006762125618);
\draw (14.52,2.94)-- (7.78,2.9);
\draw (11.045001849047301,2.919376865572981) -- (11.149198827165263,3.0549976226534703);
\draw (11.045001849047301,2.919376865572981) -- (11.150801172834738,2.785002377346529);
\draw (14.52,2.94)-- (14.56,-3.8);
\draw (14.541246268854039,-0.6399963019053985) -- (14.405625511773549,-0.5357993237874382);
\draw (14.541246268854039,-0.6399963019053985) -- (14.67562075708049,-0.5341969781179603);
\draw (14.54,-0.43) -- (14.40437924291951,-0.32580302188203975);
\draw (14.54,-0.43) -- (14.674374488226452,-0.3242006762125618);
\draw (14.56,-3.8)-- (7.82,-3.84);
\draw (11.0850018490473,-3.820623134427019) -- (11.189198827165262,-3.6850023773465295);
\draw (11.0850018490473,-3.820623134427019) -- (11.190801172834737,-3.9549976226534707);
\draw [line width=1.2pt] (7.78,2.9)-- (12.440193355815202,-3.8125804548616307);
\draw [line width=1.2pt] (12.200200399742895,2.926232643321916)-- (14.541838993422731,-0.7398703917304288);
\draw [line width=1.2pt] (7.801601613059442,-0.7398718005159677)-- (10.040159192765886,-3.8268239810518345);
\draw [line width=1.2pt] (9.800047546511935,2.911988412738943)-- (14.56,-3.8);
\draw [line width=1.2pt, color=gray] (7.7970913861569215,0.02010143255879493)-- (10.540022892765005,-3.8238574309034714);
\draw [line width=1.2pt, color=gray](8.200103897933486,2.90249319820732)-- (12.880177859174276,-3.8099692708654342);
\draw [line width=1.2pt, color=gray](10.260031345478238,2.9147182869167847)-- (14.556439294550641,-3.2000211317830822);
\draw [line width=1.2pt, color=gray] (12.560306410854693,2.9283697709843004)-- (14.538041612002852,-0.10001162248069573);
\end{tikzpicture}
\end{center}
\caption{one-sided fattening of $\omega$-limit sets in $S^2,T^2$}
\label{Fig: 1-D}
\end{figure}
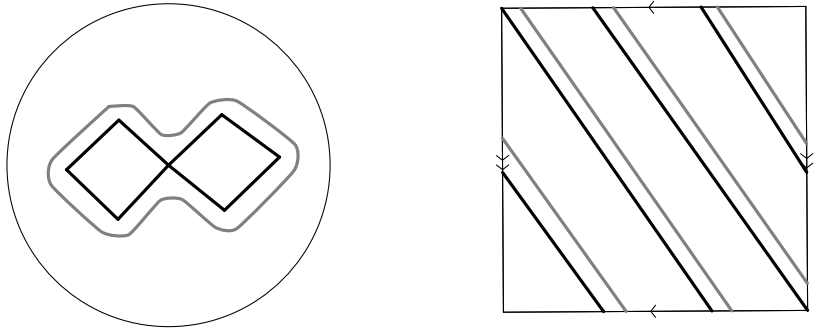

Thirdly the oriented graph $\Gamma$ corresponding to a polycycle $\omega$-limit set is not just a cactus boundary, but admits some fattening.
The definition and construction of fattening will be given precisely in \S\ref{subsec: semiribbon}.
But we give its meaning: any point in the one-sided fattened region is attracted to $\omega(x)$ as $t\to \infty$.
This explains the statement on the convergence $\Phi(t,y)\to D$ as $t\to\infty$ in the theorem.
See Fig.\ \ref{Fig: 1-D} for genus $0,1$.

In Fig.\ \ref{Fig: 1-D}, one boundary component of the fattened region defines a full Eulerian path of $\Gamma$.
Our approach to Theorem \ref{th: main: analytic version} is based on graph theory of the oriented graphs with fattening.
Hence Theorem \ref{th: main: analytic version} follows from a graph-theoretic counterpart Theorem \ref{th: number of Eulerian paths}: the pairs $(S,\Gamma)$ with fattening always come as in the statement of Theorem \ref{th: main: analytic version}.


\subsection{Motivation}\label{subsec: Motivation and further questions}

The notion of limit sets has important role to study dynamical systems. 
This expresses limiting behavior of an orbit. 
In particular the $\omega$-limit set can be used to analyze the long-term behavior of a forward orbit, which practically reconstructs the whole orbit.
This underlying idea appears in various works in a far general context in dynamical systems, e.g.\  \cite{MP, MPP, BM}.

Morales and Pacifico \cite{MP} proved that for a $\cC^{1}$-flow on a closed manifold, the set of points with Lyapunov stable $\omega$-limit set is residual in the manifold. 
For a flow on a $3$-manifold, Morales, Pacifico and Pujals \cite{MPP} clarified the structure of all $\cC^1$-robust transitive sets with singularities.
For instance, they proved that the above set is either a proper attractor or a proper repeller.
Bautista and Morales \cite{BM} found a condition equivalent to that a singular hyperbolic $\omega$-limit set in a compact manifold is a closed orbit.   

A viewpoint more relevant to ours is study of a relation between topology of a base manifold (equipped with a flow) and $\omega$-limit sets.
In general, Barwell et al.\ \cite{BGOR} found a sufficient condition under which all notions of $\omega$-limit set, weak incompressibility and internal chain transitivity become equivalent.

Our interest in this paper is a surface flow. 
For a (real) analytic flow on $\rr^{2}$, this direction has been pursued since Poincar\'e \cite{Po} 
and then Bendixson \cite{Be}.
See \cite[\S1]{JL} for the nice survey of surface flows for keeping track of notable progresses up to present.
In any case of polynomial or analytic flows, the complete topological classification of $\omega$-limit sets has been done only for $\rr^{2},S^{2}$ or the real projective plane $\rr\mathrm{P}^{2}$ (the original proof in \cite{JL} has a gap, which is reconsidered in \cite{BL}; see Remark \ref{rk: gap}).
In any cases $S^{2},\rr\mathrm{P}^{2}$, the $\omega$-limit sets turn out to be a polycycle.

There is also an approach from a general surface flow to a reasonable sufficient condition that an $\omega$-limit set is a finite union of orbits.
The following well-known theorems assure that every $\omega$-limit set is a polycycle, if a given flow is either
	\begin{itemize}
	\item
	a $\cC^{1}$-flow with only hyperbolic singularities on a compact surface and only trivial recurrent orbits  (\cite[Proposition 4.2.3]{PM}), or
	\item
	an analytic flow on a compact surface with finitely many singularities and only trivial recurrent orbits (a.k.a.\ Bendixson's theorem).
	\end{itemize}
See further explanation on these theorems in \cite[Theorems 1.4, 1.5]{JL}.
Hence, combined with our main theorem (or its generalization Theorem \ref{th: main: omega}), these yield an immediate corollary:

\begin{cor}\label{cor: Pujals de Melo Bendixson}
Let $S$ be an oriented closed surface.
For any flow on $S$ satisfying one of the above two assumptions, $(\omega(x),S)$ is given as in the statement of Theorem \ref{th: main: analytic version}, unless $\omega(x)$ is a singular point or a periodic orbit.
\end{cor}

Attractors (and repellers) are another useful limit sets in dynamical systems (see e.g.\ \cite{Co,Mi}).
Jim\'enez L\'opez and Peralta-Salas \cite{JP} gave a topological characterization of global attractors for analytic and polynomial flows on $\rr^{2}$.
For flows on $3$-dimensional manifolds, Carballo and Morales \cite{CM} studied $\omega$-limit sets in an attractor block.
They proved that for an associated vector field, any sufficiently close vector field satisfies the property that the set $\{x| \mbox{$\omega(x)$ has a singularity}\}$ is residual in the attractor block.  
We expect that there is also a topological classification of the attractors on an arbitrary closed oriented surface (equipped with an analytic flow), which will be done in the near future using our main theorem.

There is an approach to attractors which is valid in discrete dynamics.
For a closed $n$-dimensional manifold $M$ ($n \geq 3$), Grines and Zhuzhoma \cite{GZ} proved that if $M$ admits a structurally stable self-diffeomorphism with an orientable expanding attractor of codimension one, it is homeomorphic to the $n$-dimensional torus $T^{n}$.
F.\ Rodriguez Hertz and J.\ Rodriguez Hertz \cite{HH} gave a local topological and dynamical description of an expansive attractors on compact surfaces.
Moreover it is also proved that the above attractors on compact surfaces can be decomposed into invariant sets, and that the attractors are hyperbolic except possibly at a finite number of periodic points.

In this paper we use graph-theoretic approach to topological classification problem on surface flows.
Such an attempt was made in the genus $0$ case by Jim\'enez L\'opez and Llibre \cite{JL}.
We will remind their idea in the proof of Theorem \ref{th: JL}.
A new aspect of our paper is (one-sided) fattened graphs, which seem less familiar in the problems on surface flows at least for the authors.
See Remark \ref{rk: graph} for more standard notion of (two-sided) fattened graphs in the literatures.
A purely graph-theoretic question which we consider in this paper, is the full Eulerian path of an embedded graph in a Riemann surface.
See the introductory part of \S\ref{subsec: Euler}.

\vskip.3cm

\textit{Contents of the paper.}
In \S\ref{sec: surgery}, we study basic properties of surface flows and $\omega$-limit sets, in particular the changes after cutting and pasting surgery of the surface.
This will give the induction argument in the proof of the main theorem.

In \S\ref{sec: semi ribbon}, we introduce aforementioned one-sided fattening of a graph.
With this we reformulate Theorem \ref{th: main: analytic version} in term of graphs (Theorem \ref{th: number of Eulerian paths}).
The latter theorem classifies all the graphs with one-sided fattening in punctured Riemann surfaces.

In \S\ref{sec: omega}, we prove a generalized version of Theorem \ref{th: main: analytic version} for compact oriented $\cC^{1}$-surfaces with corner and $\cC^{1}$-flows with analytic sector decomposition (Definition \ref{def: analytic sector decomposition}).
In the proof we will see that $\omega(x)$ admits a semi-ribbon graph structure.

\vskip.3cm

\textit{Acknowledgement.}
The author JC is grateful to Prof.\ Mark Siggers for discussion on graph theory.


\section{Surgery of compact oriented surfaces with corner and surface
flows}\label{sec: surgery}

We review basic materials on connected sum of Riemann surfaces.
Here, connected sum is operated for multiple embedded circles.
For instance a connected sum with $S^{2}$ along two circles produces a handle (equivalently, increase genus by $1$).
We express this surgery in terms of pasting and cutting.
This gives genus induction used in the later sections.

In fact in the later sections the surgeries will be operated along invariant circles when a surface flow is given.
Some basic properties of the operation is also reviewed in \S\ref{subsec: surgery of surface flows}.


\subsection{Compact oriented surfaces with corner}

An $n$-dimensional $\cC^{1}$-\textit{manifold with corner} is a
second countable Hausdorff space with $\cC^{1}$-atlas of open subsets
$\rr^{n-k}\times \rr_{\ge0}^{k}$, $0\le k\le n$.
It is homeomorphic to a topological manifold with boundary.
An \textit{analytic manifold with corner} is defined by further imposing that the
transition functions are restrictions of analytic isomorphisms on some
open subsets in $\rr^{n}$ containing the open sets of the above atlas
respectively.

When $n=2$, we call the above space a \textit{surface with corner}.
Let $S$ be a surface with corner.
The boundary $\Bd(S)$ is homeomorphic to a disjoint union of circles.

We always assume $S$ is compact oriented.
Thus $S$ is topologically a (compact) Riemann surface with finitely many punctures (= boundary circles).
We define $g(S)$, \textit{genus} of $S$, by capping-off the boundary circles:
Let $C_{1},C_{2},...,C_{k}$ be the boundary circles of $S$.
Let $D_{1},D_{2},...,D_{k}$ be copies of a disk.
By identifying $\Bd(D_{i})$ with $C_{i}$, we obtain a closed oriented
surface (= Riemann surface), say $S_{\cl}$.
If $\Bd(S)=\emptyset$, we set $S_{\cl}$ to be $S$ itself.
Now we define $g(S):=g(S_{\cl})$.

\begin{figure}[ht]
\begin{center}
\begin{tikzpicture}[line cap=round,line join=round,>=triangle 45,x=0.7cm,y=0.7cm]
\fill[color=zzttqq,fill=zzttqq,fill opacity=0.1] (-0.76,2.32) -- (-2.32,0.78) -- (-2.334142135623731,-1.4120310216782965) -- (-0.7941421356237315,-2.972031021678297) -- (1.3978888860545653,-2.986173157302028) -- (2.9578888860545653,-1.4461731573020287) -- (2.972031021678297,0.7458578643762683) -- (1.4320310216782977,2.3058578643762684) -- cycle;
\fill[color=zzttqq,fill=zzttqq,fill opacity=0.1] (6.34,2.32) -- (4.78,0.78) -- (4.76585786437627,-1.4120310216782967) -- (6.30585786437627,-2.972031021678296) -- (8.497888886054566,-2.9861731573020274) -- (10.057888886054567,-1.4461731573020278) -- (10.072031021678296,0.745857864376269) -- (8.532031021678296,2.305857864376269) -- cycle;
\draw [color=zzttqq] (-0.76,2.32)-- (-2.32,0.78);
\draw [color=zzttqq] (-2.32,0.78)-- (-2.334142135623731,-1.4120310216782965);
\draw [color=zzttqq] (-2.334142135623731,-1.4120310216782965)-- (-0.7941421356237315,-2.972031021678297);
\draw [color=zzttqq] (-0.7941421356237315,-2.972031021678297)-- (1.3978888860545653,-2.986173157302028);
\draw [color=zzttqq] (0.4068711900599067,-2.9797794947472234) -- (0.3010024255991954,-3.1140992800045058);
\draw [color=zzttqq] (0.4068711900599067,-2.9797794947472234) -- (0.30274432483163877,-2.844104898975818);
\draw [color=zzttqq] (1.3978888860545653,-2.986173157302028)-- (2.9578888860545653,-1.4461731573020287);
\draw [color=zzttqq] (2.3273362155289408,-2.0686418192311704) -- (2.3474541252658754,-2.238480771500126);
\draw [color=zzttqq] (2.3273362155289408,-2.0686418192311704) -- (2.1577709763176314,-2.046334205033072);
\draw [color=zzttqq] (2.177888886054566,-2.216173157302028) -- (2.1980067957915,-2.386012109570984);
\draw [color=zzttqq] (2.177888886054566,-2.216173157302028) -- (2.0083236468432553,-2.1938655431039296);
\draw [color=zzttqq] (2.9578888860545653,-1.4461731573020287)-- (2.972031021678297,0.7458578643762683);
\draw [color=zzttqq] (2.972031021678297,0.7458578643762683)-- (1.4320310216782977,2.3058578643762684);
\draw [color=zzttqq] (1.4320310216782977,2.3058578643762684)-- (-0.76,2.32);
\draw [color=zzttqq] (0.23101769599465968,2.313606337445196) -- (0.33688646045537096,2.4479261227024787);
\draw [color=zzttqq] (0.23101769599465968,2.313606337445196) -- (0.3351445612229276,2.1779317416737913);
\draw (1.4320310216782977,2.3058578643762684)-- (-0.76,2.32);
\draw (-0.76,2.32)-- (-2.32,0.78);
\draw (-1.689447329474375,1.4024686619291424) -- (-1.70956523921131,1.5723076141980976);
\draw (-1.689447329474375,1.4024686619291424) -- (-1.519882090263065,1.3801610477310442);
\draw (-1.54,1.55) -- (-1.560117909736935,1.7198389522689552);
\draw (-1.54,1.55) -- (-1.37043476078869,1.5276923858019018);
\draw (-2.334142135623731,-1.4120310216782965)-- (-2.32,0.78);
\draw (-2.3263936625548043,-0.2110176959946591) -- (-2.1920738772975215,-0.31688646045536983);
\draw (-2.3263936625548043,-0.2110176959946591) -- (-2.462068258326209,-0.315144561222927);
\draw (-0.7941421356237315,-2.972031021678297)-- (-2.334142135623731,-1.4120310216782965);
\draw (-1.7116734736945887,-2.0425836922039218) -- (-1.5418345214256335,-2.022465782466987);
\draw (-1.7116734736945887,-2.0425836922039218) -- (-1.733981087892687,-2.2121489314152316);
\draw (-1.564142135623731,-2.192031021678297) -- (-1.3943031833547759,-2.1719131119413624);
\draw (-1.564142135623731,-2.192031021678297) -- (-1.5864497498218293,-2.3615962608896077);
\draw (-0.7941421356237315,-2.972031021678297)-- (1.3978888860545653,-2.986173157302028);
\draw (1.3978888860545653,-2.986173157302028)-- (2.9578888860545653,-1.4461731573020287);
\draw (2.972031021678297,0.7458578643762683)-- (2.9578888860545653,-1.4461731573020287);
\draw (2.9642825486093707,-0.45515546130736995) -- (2.829962763352088,-0.3492866968466592);
\draw (2.9642825486093707,-0.45515546130736995) -- (3.099957144380776,-0.35102859607910203);
\draw (1.4320310216782977,2.3058578643762684)-- (2.972031021678297,0.7458578643762683);
\draw (2.3495623597491546,1.376410534901894) -- (2.1797234074801985,1.356292625164959);
\draw (2.3495623597491546,1.376410534901894) -- (2.371869973947253,1.545975774113204);
\draw (2.202031021678297,1.525857864376269) -- (2.032192069409341,1.505739954639334);
\draw (2.202031021678297,1.525857864376269) -- (2.224338635876395,1.695423103587579);
\draw [color=zzttqq] (6.34,2.32)-- (4.78,0.78);
\draw [color=zzttqq] (4.78,0.78)-- (4.76585786437627,-1.4120310216782967);
\draw [color=zzttqq] (4.76585786437627,-1.4120310216782967)-- (6.30585786437627,-2.972031021678296);
\draw [color=zzttqq] (6.30585786437627,-2.972031021678296)-- (8.497888886054566,-2.9861731573020274);
\draw [color=zzttqq] (7.506871190059908,-2.979779494747222) -- (7.401002425599197,-3.114099280004505);
\draw [color=zzttqq] (7.506871190059908,-2.979779494747222) -- (7.40274432483164,-2.8441048989758166);
\draw [color=zzttqq] (8.497888886054566,-2.9861731573020274)-- (10.057888886054567,-1.4461731573020278);
\draw [color=zzttqq] (9.427336215528939,-2.0686418192311704) -- (9.447454125265875,-2.238480771500126);
\draw [color=zzttqq] (9.427336215528939,-2.0686418192311704) -- (9.257770976317628,-2.046334205033072);
\draw [color=zzttqq] (9.277888886054566,-2.216173157302028) -- (9.298006795791501,-2.386012109570984);
\draw [color=zzttqq] (9.277888886054566,-2.216173157302028) -- (9.108323646843255,-2.1938655431039296);
\draw [color=zzttqq] (10.057888886054567,-1.4461731573020278)-- (10.072031021678296,0.745857864376269);
\draw [color=zzttqq] (10.072031021678296,0.745857864376269)-- (8.532031021678296,2.305857864376269);
\draw [color=zzttqq] (8.532031021678296,2.305857864376269)-- (6.34,2.32);
\draw [color=zzttqq] (7.3310176959946585,2.313606337445196) -- (7.436886460455369,2.4479261227024787);
\draw [color=zzttqq] (7.3310176959946585,2.313606337445196) -- (7.435144561222926,2.1779317416737913);
\draw (8.532031021678296,2.305857864376269)-- (6.34,2.32);
\draw (6.34,2.32)-- (4.78,0.78);
\draw (5.410552670525625,1.4024686619291424) -- (5.39043476078869,1.5723076141980976);
\draw (5.410552670525625,1.4024686619291424) -- (5.580117909736934,1.3801610477310442);
\draw (5.56,1.55) -- (5.539882090263064,1.7198389522689552);
\draw (5.56,1.55) -- (5.72956523921131,1.5276923858019018);
\draw (4.76585786437627,-1.4120310216782967)-- (4.78,0.78);
\draw (4.773606337445195,-0.2110176959946591) -- (4.907926122702478,-0.31688646045536983);
\draw (4.773606337445195,-0.2110176959946591) -- (4.63793174167379,-0.315144561222927);
\draw (6.30585786437627,-2.972031021678296)-- (4.76585786437627,-1.4120310216782967);
\draw (5.388326526305413,-2.0425836922039218) -- (5.558165478574368,-2.0224657824669863);
\draw (5.388326526305413,-2.0425836922039218) -- (5.3660189121073145,-2.2121489314152316);
\draw (5.53585786437627,-2.192031021678297) -- (5.705696816645226,-2.171913111941361);
\draw (5.53585786437627,-2.192031021678297) -- (5.513550250178172,-2.3615962608896077);
\draw (6.30585786437627,-2.972031021678296)-- (8.497888886054566,-2.9861731573020274);
\draw (8.497888886054566,-2.9861731573020274)-- (10.057888886054567,-1.4461731573020278);
\draw (10.072031021678296,0.745857864376269)-- (10.057888886054567,-1.4461731573020278);
\draw (10.06428254860937,-0.4551554613073688) -- (9.929962763352087,-0.3492866968466581);
\draw (10.06428254860937,-0.4551554613073688) -- (10.199957144380773,-0.3510285960791009);
\draw (8.532031021678296,2.305857864376269)-- (10.072031021678296,0.745857864376269);
\draw (9.449562359749152,1.376410534901894) -- (9.279723407480196,1.356292625164959);
\draw (9.449562359749152,1.376410534901894) -- (9.47186997394725,1.545975774113204);
\draw (9.302031021678294,1.525857864376269) -- (9.132192069409339,1.505739954639334);
\draw (9.302031021678294,1.525857864376269) -- (9.324338635876392,1.695423103587579);
\draw [color=black] (5.9,-0.32) circle (0.8836288813749809);
\draw [color=black] (8.74,-0.32) circle (0.8836288813749809);
\fill[color=white,fill=white] (5.9,-0.32) circle (0.8736288813749809);
\fill[color=white,fill=white] (8.74,-0.32) circle (0.8736288813749809);
\begin{scriptsize}
\draw[color=black] (0.38,2.8) node {$\alpha_1$};
\draw[color=black] (-1.72,1.94) node {$\beta_1$};
\draw[color=black] (-1.76,-0.14) node {$\alpha_1\inv$};
\draw[color=black] (-1.08,-1.8) node {$\beta_1\inv$};
\draw[color=black] (0.34,-3.32) node {$\alpha_2$};
\draw[color=black] (2.66,-2.28) node {$\beta_2$};
\draw[color=black] (2.5,-0.18) node {$\alpha_2\inv$};
\draw[color=black] (1.72,1.10) node {$\beta_2\inv$};
\end{scriptsize}
\end{tikzpicture}
\end{center}
\caption{planar diagrams of $S$ with $g(S_{\cl})=2$; there are no (resp.\ two) punctures in the left (resp.\ right)}
\label{Fig: 2-A}
\end{figure}
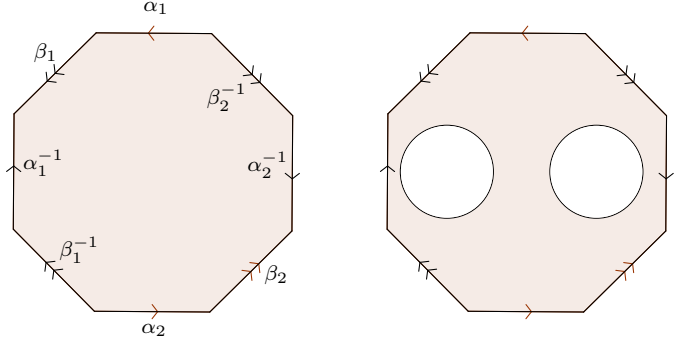

It is well-known that $S_{\cl}$ has a planar diagram given by
$4g(S_{\cl})$-gon if $S_{\cl}$ is connected and $g(S_{\cl})\ge1$.
The edges of the $4g(S_{\cl})$-gon are labelled by
$\alpha_{i},\beta_{i},\alpha_{i}\inv,\beta_{i}\inv$, $1\le i\le
g(S_{\cl})$ counterclockwise (Fig.\ \ref{Fig: 2-A}).

$S_{\cl}$ is obtained by identifying $\alpha_{i},\alpha_{i}\inv$ with
opposite orientation and similarly $\beta_{i},\beta_{i}\inv$.
If $S_{\cl}$ is connected and $g(S_{\cl})\ge1$, a planar diagram of $S$
is given by taking finitely many punctures in the interior of the above
planar diagram of $S_{\cl}$.
Even when $g(S_{\cl})=0$, if $\Bd(S)\neq\emptyset$, there is a planar
diagram of $S$.
It is a disk with $k-1$ punctures where $k$ denotes the number of punctures of $S$.

Let $H_{i}(S,\zz)$ be the $i^{\mathrm{th}}$ singular homology group of $S$ over $\zz$.
Suppose $S$ is connected.
Then $H_{0}(S,\zz)\cong\zz$ by the map sending the class of a point to $1\in \zz$.
If $S=S_{\cl}$, $H_{2}(S,\zz)\cong \zz$ by the map sending the fundamental class $[S]$ to $1\in\zz$.
In this case $H_{1}(S,\zz)=\zz\langle [\alpha_{i}],[\beta_{i}]\rangle_{i=1,2,...,g(S)}$ the $\zz$-module freely generated by the classes $[\alpha_{i}],[\beta_{i}]$.

Assume $S\neq S_{\cl}$, i.e., there are the punctures $C_{1},C_{2},...,C_{k}$ on $S$. 
We give any orientation on $C_{i}$ for each $i$.
Then we have 
	\begin{equation}\label{eq: bouquet}
	\begin{aligned}
	&
	H_{2}(S,\zz)=0,
	\\
	&
	H_{1}(S,\zz)=\zz\langle [\alpha_{i}],[\beta_{i}],[C_{j}]\rangle_{i=1,2,...,g(S),\ j=1,2,...,k-1}.
 	\end{aligned}
	\end{equation}
One needs care on the indices $1\le j\le k-1$ in the above.
These identifications follow from the fact that $S$ is deformation-retract to a bouquet of $2g(S)+k-1$ circles.


\subsection{Surgery}

We define a surgery of a compact oriented surface $S$ with corner by
cutting along a simple closed curve in the interior $\Int(S)$.
Hereafter a simple closed curve is always assumed to be piecewise smooth.

\begin{lem}\label{lem: neck}
Let $C$ be a simple closed curve in $\Int(S)$.
Then there exists an embedding $C\times(-\ep,\ep)\to \Int(S)$ for
$0<\ep\ll 1$.
\end{lem}

\proof
This is immediate using the planar diagram.
We give an alternative intrinsic proof.
By smoothing, $C$ can be assumed to be a smooth curve.
Let $N_{C}$ be the normal bundle of $C$ in $S$.
Using the orietations of $S$ and $C$, $N_{C}$ has the induced orientation.
So $N_{C}\stackrel{\mathrm{diffeo}}\approx C\times\rr^{1}$.
Using the local embedding of $N_{C}$ into $S$, we are done.
\qed\vskip.3cm

\begin{figure}[ht]
\begin{center}
\begin{tikzpicture}[line cap=round,line join=round,>=triangle 45,x=.6cm,y=.6cm]
\draw plot[smooth, tension=.7pt] coordinates
{(-3.2,-0.4) (-2.,0.)  (-1.2,0.8)  (0.,1.2)};
\draw plot[smooth, tension=.7pt] coordinates
{(1.6,1.2) (2.8,0.8) (3.6,0.)  (4.76,-0.46)};
\draw plot[smooth, tension=.7pt] coordinates
{ (-3.2,-1.6) (-2.,-2.) (-1.2,-2.8) (0.,-3.2)  (1.6,-3.2) (2.8,-2.8)  (3.6,-2.) (4.8,-1.6)};
\draw plot[smooth, tension=.7pt] coordinates
{ (0.,-0.4) (-0.4,-0.8) (-0.4,-1.2)  (0.,-1.6) (0.8,-1.84) (1.6,-1.6) (2.,-1.2) (2.,-0.8) (1.6,-0.4)};
\draw [->] (5.2,0.4) -- (6.8,0.4);
\draw [rotate around={90.:(0.,0.4)}] (0.,0.4) ellipse (0.7717512552629896 and 0.66);
\draw [rotate around={90.:(1.6,0.4)}] (1.6,0.4) ellipse (0.8 and 0.6928203230275507);
\draw plot[smooth, tension=.7pt] coordinates
{(7.2,-0.4) (8.4,0.) (9.2,0.8) (10.4,1.2) (12.,1.2) (13.2,0.8) (14.,0.) (15.16,-0.46)};
\draw plot[smooth, tension=.7pt] coordinates
{ (7.2,-1.6) (8.4,-2.) (9.2,-2.8) (10.4,-3.2) (12.,-3.2) (13.2,-2.8) (14.,-2.) (15.2,-1.6)};
\draw plot[smooth, tension=.7pt] coordinates
{ (10.4,-0.4) (10.,-0.8) (10.,-1.2) (10.4,-1.6) (11.2,-1.84) (12.,-1.6) (12.4,-1.2) (12.4,-0.8) (12.,-0.4) (11.16,-0.16) (10.4,-0.4)};
\begin{scriptsize}
\draw[color=black] (0.,2.) node {$C_1$};
\draw[color=black] (1.6,2.) node {$C_2$};
\draw[color=black] (6.,1) node {$\pi$};
\end{scriptsize}
\end{tikzpicture}
\end{center}
\caption{cutting-pasting surgery}
\label{Fig: 2-B}
\end{figure}
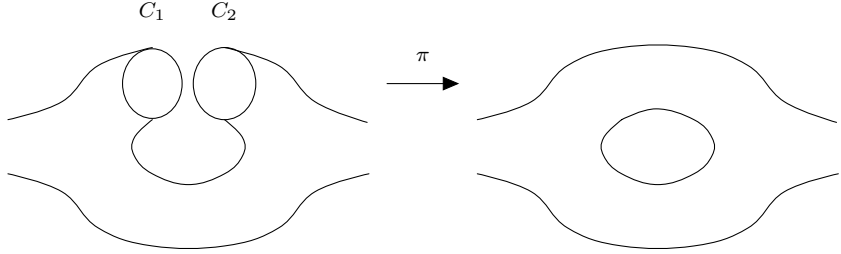

For a simple closed curve $C$ as above,  we define $\tS$ as a disjoint union
	$$
	\tS:=(S\setminus C)\sqcup C_{1}\sqcup C_{2}
	$$
where $C_{1},C_{2}$ are copies of $C$.
$\tS$ is naturally a (compact oriented) surface with corner due to Lemma
\ref{lem: neck}: the atlas of $\tS$ consists of three kinds of open
sets, (i) $U\cap (S\setminus C)$, (ii) $U\cap (C\times (-\ep,0])$, (iii)
$U\cap (C\times[0,\ep))$ where $U$ is an element of the atlas of $S$.
The process from $S$ to $\tS$ is \textit{cutting surgery} and the reverse process is \textit{pasting surgery}.
The obvious map $\pi\colon \tS\to S$ corresponds to these surgeries depicted as Fig.\ \ref{Fig: 2-B}.

We notice that both $C_{1},C_{2}$ are new punctures  of $\tS$.
So we can take a planar diagram of $\tS$ containing $C_{1},C_{2}$ in its
interior.
But we need care as $\tS$ can be disconnected even if $S$ is connected.

\begin{lem}\label{lem: conn comp}
The number of connected components of $\tS$ equals that of $S$ or
increase by $1$.
\end{lem}

\proof
We may assume $S$ is connected because the components of $S$ not
containing $C$, do not change by the above cut surgery.
Every component of $\tS$ contains $C_{1}$ or $C_{2}$ (or both) as
punctures.
Two or more components do not contain one of $C_{i}$ since otherwise
$C_{i}$ cannot be a puncture.
Therefore there are only two possibilities: one component contains both
or there are two components containing $C_{1},C_{2}$ respectively.
\qed\vskip.3cm
Note that in the latter case of the above lemma, $\tS$ is a connected sum along one puncture.

In the rest of this subsection we compute the genus of $\tS$.
This will serve as genus induction for the proofs in foregoing claims.
We may assume $S$ is connected.

Let us consider first the case $\tS$ is connected.
Since $C_{1},C_{2}$ are punctures of $\tS$, there is a planar diagram as
in Fig.\ \ref{Fig: 2-C}.
\begin{figure}[ht]
\begin{center}
\begin{tikzpicture}[line cap=round,line join=round,>=triangle 45,x=0.6cm,y=0.6cm]
\fill[color=zzttqq,fill=zzttqq,fill opacity=0.1] (2.,-4.8) -- (6.4,-4.8) -- (9.511269837220809,-1.688730162779191) -- (9.511269837220809,2.711269837220809) -- (6.4,5.822539674441619) -- (2.,5.82253967444162) -- (-1.1112698372208092,2.7112698372208106) -- (-1.11126983722081,-1.6887301627791893) -- cycle;
\draw [color=zzttqq] (2.,-4.8)-- (6.4,-4.8);
\draw [color=zzttqq] (6.4,-4.8)-- (9.511269837220809,-1.688730162779191);
\draw [color=zzttqq] (9.511269837220809,-1.688730162779191)-- (9.511269837220809,2.711269837220809);
\draw [color=zzttqq] (9.511269837220809,2.711269837220809)-- (6.4,5.822539674441619);
\draw [color=zzttqq] (6.4,5.822539674441619)-- (2.,5.82253967444162);
\draw [color=zzttqq] (4.095,5.822539674441619) -- (4.2,5.957539674441619);
\draw [color=zzttqq] (4.095,5.822539674441619) -- (4.2,5.687539674441618);
\draw [color=zzttqq] (2.,5.82253967444162)-- (-1.1112698372208092,2.7112698372208106);
\draw [color=zzttqq] (-1.1112698372208092,2.7112698372208106)-- (-1.11126983722081,-1.6887301627791893);
\draw [color=zzttqq] (-1.11126983722081,-1.6887301627791893)-- (2.,-4.8);
\draw(2.,2.) circle (1.2);
\draw(6.4,2.) circle (1.2);
\draw(4.,-2.) circle (1.2);
\fill[color=white,fill=white] (2.,2.) circle (1.19);
\fill[color=white,fill=white ] (6.4,2.) circle (1.19);
\fill[color=white,fill=white ] (4.,-2.) circle (1.19);
\draw (6.4,5.822539674441619)-- (2.,5.82253967444162);
\draw (2.,5.82253967444162)-- (-1.1112698372208092,2.7112698372208106);
\draw (0.2958726573404209,4.11841233178204) -- (0.2746594539048243,4.2881179592668115);
\draw (0.2958726573404209,4.11841233178204) -- (0.4655782848251926,4.0971991283464435);
\draw (0.44436508138959596,4.266904755831215) -- (0.4231518779539994,4.436610383315987);
\draw (0.44436508138959596,4.266904755831215) -- (0.6140707088743677,4.245691552395618);
\draw (-1.11126983722081,-1.6887301627791893)-- (-1.1112698372208092,2.7112698372208106);
\draw (-1.1112698372208099,0.6162698372208104) -- (-0.9762698372208098,0.5112698372208103);
\draw (-1.1112698372208099,0.6162698372208104) -- (-1.2462698372208099,0.5112698372208103);
\draw (2.,-4.8)-- (-1.11126983722081,-1.6887301627791893);
\draw (0.2958726573404198,-3.0958726573404194) -- (0.4655782848251914,-3.0746594539048226);
\draw (0.2958726573404198,-3.0958726573404194) -- (0.2746594539048232,-3.265578284825191);
\draw (0.44436508138959485,-3.2443650813895943) -- (0.6140707088743664,-3.2231518779539976);
\draw (0.44436508138959485,-3.2443650813895943) -- (0.4231518779539982,-3.414070708874366);
\draw (2.,-4.8)-- (6.4,-4.8);
\draw (4.305,-4.8) -- (4.2,-4.935);
\draw (4.305,-4.8) -- (4.2,-4.665);
\draw (6.4,-4.8)-- (9.511269837220809,-1.688730162779191);
\draw (8.104127342659579,-3.0958726573404203) -- (8.125340546095174,-3.265578284825192);
\draw (8.104127342659579,-3.0958726573404203) -- (7.934421715174808,-3.074659453904824);
\draw (7.955634918610404,-3.2443650813895957) -- (7.9768481220459995,-3.414070708874367);
\draw (7.955634918610404,-3.2443650813895957) -- (7.785929291125633,-3.223151877953999);
\draw (9.511269837220809,2.711269837220809)-- (9.511269837220809,-1.688730162779191);
\draw (9.511269837220809,0.4062698372208092) -- (9.376269837220809,0.5112698372208092);
\draw (9.511269837220809,0.4062698372208092) -- (9.64626983722081,0.5112698372208092);
\draw (6.4,5.822539674441619)-- (9.511269837220809,2.711269837220809);
\draw (8.104127342659579,4.1184123317820385) -- (7.934421715174808,4.097199128346443);
\draw (8.104127342659579,4.1184123317820385) -- (8.125340546095174,4.28811795926681);
\draw (7.955634918610404,4.266904755831214) -- (7.785929291125633,4.245691552395617);
\draw (7.955634918610404,4.266904755831214) -- (7.9768481220459995,4.436610383315985);
\begin{scriptsize}
\draw[color=black] (2.14,3.56) node {$C_1$};
\draw[color=black] (6.54,3.56) node {$C_2$};
\end{scriptsize}
\end{tikzpicture}
\end{center}
\caption{}
\label{Fig: 2-C}
\end{figure}

In Fig.\ \ref{Fig: 2-C} the readers should be aware of that orientations of $C_{1},C_{2}$ are given opposite.
The identification along the orientation (pasting surgery), we
recover $C$.
The pasting surgery is nothing but taking the image of $\pi\colon \tS\to S$.

\begin{lem}\label{lem: genus: connected case}
If $\tS$ is connected,  $g(\tS)=g(S)-1$.
\end{lem}

\proof
From the planar diagram of Fig.\ \ref{Fig: 2-C}, we see the pasting surgery increases genus by $1$.
\qed\vskip.3cm

Next we consider the case $\tS$ is disconnected.
We denote the two connected components by $\tS_{1},\tS_{2}$ with
$C_{i}\subset \tS_{i}$.
In this case a similar pasting surgery is also nothing but
identification of two planar diagrams along $C_{1},C_{2}$.

\begin{lem}\label{lem: genus: disconnected case}
$g(\tS)=g(\tS_{1})+g(\tS_{2})$.
Moreover if $[C]\neq0$ in the singular homology group
$H_{1}(S_{\cl},\zz)$, both $g(\tS_{1}),g(\tS_{2})$ are strictly less
than $g(\tS)$.
\end{lem}

In the statement the homology class $[C]$ is defined after fixing any orientation  
of $C$.
The assumption $[C]\neq0$ does not depend on the choice of orientation.
\vskip.3cm

\textit{Proof of Lemma \ref{lem: genus: disconnected case}.}
To prove the first assertion we add one puncture to $\tS_{1},\tS_{2}$ respectively. 
Note that the genus do not depend on the number of punctures.
By the same symbols $\tS_{1},\tS_{2}$, we denote the Riemann surfaces with the additional punctures.
By \eqref{eq: bouquet}, $g(\tS)$ is the half of the number of the circles in the bouquets which are not punctures (\textit{non-puncture circles} for short).
On the other hand $\tS$ is a connected sum along $C_{1},C_{2}$.
Thus the bouquet of $\tS$ is a union of the bouquets of $\tS_{1},\tS_{2}$ with the two puncture circles $C_{1},C_{2}$ being identified. 
It is obvious that the total number of non-puncture circles does not change under this identification.
This proves  the identity $g(\tS)=g(\tS_{1})+g(\tS_{2})$.

Let us prove the second assertion.
Now we suppose the contrary, say $g(\tS_{1})=0$.
Then $C_{1}$ is contractible to a point in the surface obtained from $\tS_{1}$ by capping off all the punctures other than $C_{1}$.
Thus $C$ is also contractible to a point in $S_{\cl}$.
This contradicts $[C]\neq0$.
\qed\vskip.3cm


\subsection{Surgery of surface flows}\label{subsec: surgery of surface flows}

Let $\Phi$ be a flow on a punctured Riemann surface $S$.
Let us assume $C\subset S$ be a $\Phi$-invariant simple closed curve.
Let $\tS$ be the cutting surgery of $S$ along $C$.

\begin{lem}\label{lem: lift of flow}
There is a unique flow $\tPhi$ on $\tS$ which fits in the commutative diagram:
	$$
	\xymatrix{
	\rr\times\tS\ar[r]^-{\tPhi}\ar[d]_{\mathrm{Id}_{\rr}\times \pi} & \tS\ar[d]^{\pi}
	\\
	\rr\times S\ar[r]^-{\Phi} & S
	}
	$$
Moreover $\tPhi$ preserves given structures of $(S,\Phi)$ (e.g.\ $\cC^{1}$, smoothness,
analyticity).
\end{lem}

\proof
We define $\tPhi(t,x):=\Phi(t,x)$ by identifying $C=C_{1}=C_{2}$.
The structures are preserved due to the embedding $N_{C}\to S$ also
preserves the structure.
\qed\vskip.3cm

We will obtain comparison result on $\omega$-limit sets for $(S,\Phi)$ and
$(\tS,\tPhi)$.
Recall the \textit{$\omega$-limit set} $\omega(x)$ of $x$ is defined by
    $$
\omega(x):= \left\{ z \in M \middle| z=\lim_{n\to\infty}
\Phi(t_{n},x)
\mbox{ \ for some sequence $t_{n}\to \infty$ as $n\to \infty$}\right\}.
    $$
The $\alpha$-limit set $\alpha(x)$ is similarly defined replacing $t_{n}\to \infty$ into $t_{n}\to -\infty$ in the above.

\begin{lem}\label{lem: lift and omega}
Let $x\in S\setminus C$ and $\tx:=\pi\inv(x)$.
Then $\pi(\omega(\tx))=\omega(x)$.
\end{lem}

\proof
Since $x\in S\setminus C$, we have $\pi(\tPhi(t,\tx))=\Phi(t,x)$ for $t\in \rr$.
Thus $\pi(\omega(\tx))\subset\omega(x)$.

We check the opposite inclusion. 
Let $y\in \omega(x)$.
If $y\notin C$, the opposite inclusion is clear since $\pi$ is a homeomorphism near $y$.
Otherwise it is also clear since $\pi$ is a two-one map over $C$.
\qed\vskip.3cm

This lemma does \textit{not} imply $\omega(\tx)=\pi\inv(\omega(x))$.
This phenomenon happens when $\omega(x)$ contains $C$.
Later we completely analyze the $\omega$-limit sets in this case during the proof of Theorem \ref{th: main: omega}.


\section{Semi-ribbon graph and Eulerian path}\label{sec: semi ribbon}

In this section we develop a graph-theoretic counterpart of $\omega$-limit sets.
In \S\ref{subsec: semiribbon} we introduce semi-ribbon graph (= one-sided fattening) of a graph $\Gamma$ embedded in a surface $S$. 
In this section our discussion is always topological so that $S$ is a topological surface with boundary.

The main purpose in this section is to classify pairs $(S,\Gamma)$ allowing semi-ribbon graph (Theorem \ref{th: number of Eulerian paths}).
The proof of the theorem in \S\ref{subsec: Euler} uses genus induction via cutting surgery along a simple closed path of $\Gamma$ in the previous section.
The initial induction hypothesis is the case when $g(S)=0$ (Lemma \ref{lem: geneus 0}).
It turns out that the one-sided fattening always comes from an Eulerian path on $\Gamma$.


\subsection{Semi-ribbon graph}\label{subsec: semiribbon}

An \textit{oriented graph} $\Gamma=(V,E)$ is a pair of nonempty sets with
two functions $h,t\colon E\to V$.
The elements of $V,E$ are \textit{vertices} and \textit{edges}
respectively.
For each $e\in E$, the vertices $h(e),t(e)$ are \textit{head} and
\textit{tail} vertices of $e$.
If $h(e)=t(e)$, $e$ is a \textit{loop}.
\textit{Valency} at a vertex $v$ of $\Gamma$ is the sum
$\#h\inv(v)+\#t\inv(v)$, the number of edges having $v$ as head or tail
(counted twice when loop).
We always assume that $\Gamma$ has finite valency at each vertex.
In paper a graph always means an oriented graph for short.
But sometimes we use `oriented graph' to emphasize orientation.

We frequently regard $\Gamma$ as a topological space (with orientation) defined as follows.
First we give discrete topology on $E$.
Let $[0,1]_{e}:=[0,1]\times\{e\}\subset \rr^{1}\times E$.
Let $t_{e}:=(t,e)\in [0,1]_{e}$, $0\le t\le 1$.
We define an equivalence relation on the disjoint union $\bigsqcup_{e\in E}[0,1]_{e}$ such that the only nontrivial relations are
	$$
	1_{e}\sim 0_{e'}\ \mbox{if and only if}\ h(e)=t(e') 
	$$
First we consider the line segments $l_{e}\approx[0,1]$, $e\in E$.
Now an element of $\Gamma$ makes sense as $\Gamma$ is a set.
A \textit{smooth point} of $\Gamma$ means an element of $\Gamma$ which does not lie in $V$.

Let $S$ be a surface with boundary not necessarily oriented.

\begin{defn}\label{def: hyperbolic graph}
An oriented graph $\Gamma$ embedded in $S$ is \textit{hyperbolic} if
each vertex of $\Gamma$ has an open neighborhood consisting of
alternating orientations as in Fig.\ \ref{Fig: 3-A}.
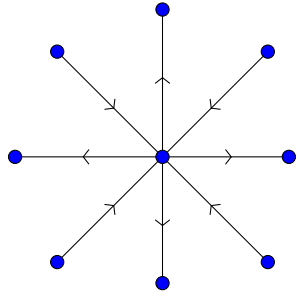
\begin{figure}[ht]
\begin{center}
\begin{tikzpicture}[line cap=round,line join=round,>=triangle 45,x=0.7cm,y=.7cm]
\draw (5.6,0.4)-- (5.6,3.2);
\draw (5.6,1.905) -- (5.735,1.8);
\draw (5.6,1.905) -- (5.465,1.8);
\draw (3.6,2.4)-- (5.6,0.4);
\draw (4.674246212024587,1.3257537879754124) -- (4.504540584539816,1.304540584539816);
\draw (4.674246212024587,1.3257537879754124) -- (4.695459415460184,1.4954594154601841);
\draw (5.6,0.4)-- (2.8,0.4);
\draw (4.095,0.4) -- (4.2,0.535);
\draw (4.095,0.4) -- (4.2,0.265);
\draw (5.6,0.4)-- (5.6,-2.);
\draw (5.6,-0.905) -- (5.465,-0.8);
\draw (5.6,-0.905) -- (5.735,-0.8);
\draw (5.6,0.4)-- (8.,0.4);
\draw (6.905,0.4) -- (6.8,0.265);
\draw (6.905,0.4) -- (6.8,0.535);
\draw (7.6,2.4)-- (5.6,0.4);
\draw (6.525753787975414,1.3257537879754124) -- (6.504540584539818,1.4954594154601841);
\draw (6.525753787975414,1.3257537879754124) -- (6.695459415460184,1.304540584539816);
\draw (3.6,-1.6)-- (5.6,0.4);
\draw (4.674246212024587,-0.5257537879754125) -- (4.695459415460184,-0.6954594154601841);
\draw (4.674246212024587,-0.5257537879754125) -- (4.504540584539816,-0.5045405845398159);
\draw (7.6,-1.6)-- (5.6,0.4);
\draw (6.525753787975414,-0.5257537879754125) -- (6.695459415460184,-0.5045405845398159);
\draw (6.525753787975414,-0.5257537879754125) -- (6.504540584539818,-0.6954594154601841);
\begin{scriptsize}
\draw [fill=qqqqff] (5.6,0.4) circle (2.5pt);
\draw [fill=qqqqff] (5.6,3.2) circle (2.5pt);
\draw [fill=qqqqff] (3.6,2.4) circle (2.5pt);
\draw [fill=qqqqff] (2.8,0.4) circle (2.5pt);
\draw [fill=qqqqff] (5.6,-2.) circle (2.5pt);
\draw [fill=qqqqff] (8.,0.4) circle (2.5pt);
\draw [fill=qqqqff] (7.6,2.4) circle (2.5pt);
\draw [fill=qqqqff] (3.6,-1.6) circle (2.5pt);
\draw [fill=qqqqff] (7.6,-1.6) circle (2.5pt);
\end{scriptsize}
\end{tikzpicture}
\end{center}
\caption{hyperbolic at a vertex}
\label{Fig: 3-A}
\end{figure}
\end{defn}

This definition makes sense even when an edge lies in the boundary of $S$.
Note any hyperbolic graph has even valency at any vertex.

Now we introduce semi-ribbon graph of $\Gamma$ via the aforementioned fattening $\Gamma$ in \S\ref{subsec: main theorem}.
For the purpose we define first some orientation rule on the sectors given by $(\Gamma,S)$ where $S$ is a (oriented) Riemann surface with boundary and $\Gamma$ is embedded in $S$. 
Let $B_{0}$ be the open unit disk in $\rr^{2}$ and $B_{0}^{I}$ be the subset in $B_{0}$ of $(x_{1},x_{2})$ with $x_{2}\in I$.
E.g., $B_{0}^{\ge0}$ is the upper half open unit disk.
For any point $x\in\Gamma$, there is an open neighborhood $U_{x}$ of $x$ in $S$ with an orientation-preserving based homeomorphism $\phi_{x}\colon (U_{x},x)\approx (B_{0},(0,0))$ such that the image $\phi_{x}(U_{x}\cap\Gamma)$ becomes a sector decomposition of $B_{0}$ centered at the origin $(0,0)\in \rr^{2}$.
We also say that  $U_{x}\cap \Gamma$ is a sector decomposition of $U_{x}$ centered at $x$.
For instance the graph in Fig.\ \ref{Fig: 3-A} has the eight sectors at the central vertex.
If $x$ is a smooth point of $\Gamma$, $U_{x}$ is decomposed into two sectors.

A sector has two line segments called \textit{separatices}.
The two separatices coincide if and only if the angle between them is $360^{\circ}$.
This case occurs only when $x$ is a vertex of $\Gamma$ with valency $1$.

For any fixed sector $T$ in $U_{x}$ with distinct separatices, we consider a homeomorphism $\psi$ from $(B_{0},(0,0))$ to itself such that $\psi(\phi_{x}(T))=B_{0}^{\ge0}$.
By replacing $\phi_{x}$ to the composite $\psi\phi_{x}$, we may assume $\phi_{x}$ satisfies $\phi_{x}(T)=B_{0}^{\ge0}$.
It is clear that for a given sector $T$ centered at $x$, the existence of such a homeomorphism $\phi_{x}$ assures that $T$ has two distinct separatrices.

\begin{defn}\label{def: oriented sector}
A sector $T$ in $U_{x}$ is \textit{oriented} if the two separatices are distinct and the arrows associated to them are heading from left to right via $\phi_{x}$. 
See Fig.\ \ref{Fig: 3-a}.

\begin{figure}[ht]
\begin{center}
\begin{tikzpicture}[line cap=round,line join=round,>=triangle 45,x=.6cm,y=.6cm]
\draw (1.2,0.4)-- (4.4,0.4);
\draw (2.905,0.4) -- (2.8,0.265);
\draw (2.905,0.4) -- (2.8,0.535);
\draw (4.4,0.4)-- (7.6,0.4);
\draw (6.105,0.4) -- (6.,0.265);
\draw (6.105,0.4) -- (6.,0.535);
\draw [dash pattern=on 5pt off 5pt] (4.4,0.4) circle (3.2);
\begin{scriptsize}
\draw [fill=qqqqff] (1.2,0.4) circle (2.5pt);
\draw [fill=qqqqff] (4.4,0.4) circle (2.5pt);
\draw [fill=qqqqff] (7.6,0.4) circle (2.5pt);
\end{scriptsize}
\end{tikzpicture}
\end{center}
\caption{}
\label{Fig: 3-a}
\end{figure}
\end{defn}

Note that the orientation of a sector is well-defined.
In other words for any choice of $\phi_{x}$, the orientation of the boundary $\Bd(T)$ induced by $S$ coincides with the one of the upper half unit disk.
This is due to the orientation-preserving assumption of $\phi_{x}$.

We define a pre-version of semi-ribbon graph first.

\begin{defn}\label{def: pre semi ribbon}
Let $(S,\Gamma)$ be a pair of an oriented graph $\Gamma$ embedded in a Riemann surface $S$ with boundary.
A closed subset in $S$ is a \textit{pre-semi-ribbon graph} $\Gamma_{\psr}$ of $\Gamma$ if there is a set of open neighborhoods $U_{x}$ of $x\in \Gamma$ in the closed Riemann surface $S_{\cl}$ covering $\Gamma_{\psr}$ such that for each oriented sector $T$ centered at $x$, there exists an orientation-preserving homeomorphism $\phi_{x}\colon (U_{x},T)\approx (B_{0},B_{0}^{\ge0})$ satisfying 
    \begin{equation}   \label{eq: pre semi ribbon}
	\phi_{x}(\Gamma_{\psr}\cap T)=B_{0}^{[0,1/2]},\quad 
	\phi_{x}(\Bd(\Gamma_{\psr})\cap T)=B_{0}^{0}\sqcup B_{0}^{1/2}
    	\end{equation}
and for each non-oriented sector $T$, $\Gamma_{\psr}\cap \Int(T)=\emptyset$.
\end{defn}

The intersection of $\Gamma_{\psr}$ and an oriented sector is depicted as the shaded region in Fig.\ \ref{Fig: 3-aa}.
\begin{figure}[ht]
\begin{center}
\begin{tikzpicture}[line cap=round,line join=round,>=triangle 45,x=.6cm,y=.6cm]
\draw (1.2,0.4)-- (4.4,0.4);
\draw (4.4,0.4)-- (7.6,0.4);
\draw [dash pattern=on 5pt off 5pt] (4.4,0.4) circle (3.2);
\draw (1.6138028826977604,1.9738823410714934)-- (7.1929587447548755,1.9618519296326609);
\fill[color=zzttqq,fill=zzttqq,fill opacity=0.1] (4.4,0.4) -- (1.6138028826977604,1.9738823410714934) -- (7.1929587447548755,1.9618519296326609) -- cycle;
\fill [shift={(4.4,0.4)},color=zzttqq,fill=zzttqq,fill opacity=0.1]  (0,0) --  plot[domain=2.6273929128525024:3.141592653589793,variable=\t]({1.*3.2*cos(\t r)+0.*3.2*sin(\t r)},{0.*3.2*cos(\t r)+1.*3.2*sin(\t r)}) -- cycle ;
\fill [shift={(4.4,0.4)},color=zzttqq,fill=zzttqq,fill opacity=0.1]  (0,0) --  plot[domain=0.:0.5098871195844,variable=\t]({1.*3.2*cos(\t r)+0.*3.2*sin(\t r)},{0.*3.2*cos(\t r)+1.*3.2*sin(\t r)}) -- cycle ;
\end{tikzpicture}
\end{center}
\caption{}
\label{Fig: 3-aa}
\end{figure}
Note that in the above definition, $B_{0}^{0}\, (=[0,1]\times\{0\})$ is the union of two separatrices. 
Thus the union of $\phi_{x}\inv(B_{0}^{0})$ over all $(x,\phi_{x})$ is nothing else but the given graph $\Gamma$.
This is a subset of the boundary $\Bd(\Gamma_{\psr})$.
From \eqref{eq: pre semi ribbon}, the complement of $\Gamma$ in $\Bd(\Gamma_{\psr})$ is the union of $\phi_{x}\inv(B_{0}^{1/2})$ over all $(x,\phi_{x})$.
We denote  the first subset $\Gamma$ and its complement by $\Bd_{\sing}(\Gamma_{\psr})$ and $\Bd_{\sm}(\Gamma_{\psr})$respectively.
We call them \textit{singular boundary} and \textit{smooth boundary}respectively.
Note that the smooth boundary $\Bd_{\sm}(\Gamma_{\psr})$ is a $1$-dimensional manifold.

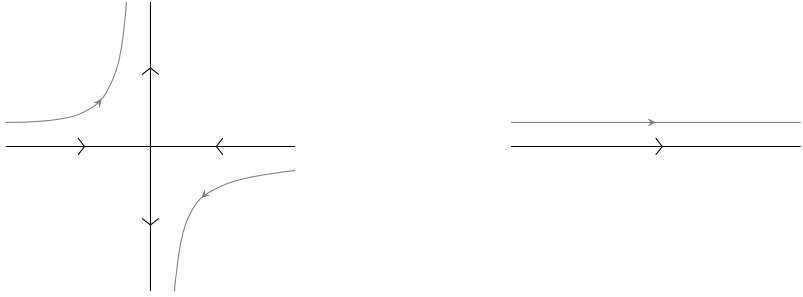
\begin{figure}[ht]
\begin{center}
\begin{tikzpicture}[line cap=round,line join=round,>=triangle 45,x=0.8cm,y=0.8cm]
\draw (10.,0.4)-- (14.8,0.4);
\draw (12.505,0.4) -- (12.4,0.265);
\draw (12.505,0.4) -- (12.4,0.535);
\draw (1.6,0.4)-- (4.,0.4);
\draw (2.905,0.4) -- (2.8,0.265);
\draw (2.905,0.4) -- (2.8,0.535);
\draw (4.,0.4)-- (4.,2.8);
\draw (4.,1.705) -- (4.135,1.6);
\draw (4.,1.705) -- (3.865,1.6);
\draw (6.4,0.4)-- (4.,0.4);
\draw (5.095,0.4) -- (5.2,0.535);
\draw (5.095,0.4) -- (5.2,0.265);
\draw (4.,0.4)-- (4.,-2.);
\draw (4.,-0.905) -- (3.865,-0.8);
\draw (4.,-0.905) -- (4.135,-0.8);
\draw[color=gray,arrow data={.5}{stealth}] plot[color=gray, smooth, tension=.7pt] coordinates
{(1.6,0.8)  (2.82,0.94) (3.4,1.58) (3.6,2.8)};
\draw[color=gray,arrow data={.5}{stealth}] plot[color=gray, smooth, tension=.7pt] coordinates
{(6.4,0.) (5.22,-0.24) (4.62,-0.8) (4.4,-2.)};
\draw[color=gray,arrow data={.5}{stealth}] (10.,0.8)-- (14.8,0.8);
\end{tikzpicture}
\end{center}
\caption{The boundaries $\Bd_{\sm}(\Gamma_{\psr})$ with orientation near a vertex with valency $4$ and a smooth point}
\label{Fig: 3-B}
\end{figure}

The arrows on the separatrices on an oriented sector induce an orientation of $\Bd_{\sm}(\Gamma_{\psr})$.
This orientation is opposite to the natural orientation of the boundary $\Bd(\Gamma_{\psr})$ induced from the one of $S$ (see the gray-colored curves in Fig.\ \ref{Fig: 3-B}).

\begin{rk}\label{rk: graph}
(1)
$\Gamma_{\psr}$ can be seen a one-sided fattening of $\Gamma$.
A similar but two-sided fattening, usually called ribbon graph, is not our interest in this paper.
This notion appears in various branches in mathematics and theoretic
physics.

(2)
A ribbon graph has the two smooth boundaries as it is a two-sided fattening.
Usually the orientation of the smooth boundaries is given opposite to our choice of the orientation of $\Bd_{\sm}(\Gamma_{\psr})$.
I.e., it obeys the natural orientation from $S$.

(3)
Our choice of orientation here is nothing but flow direction, as we will see in the proof of Theorem \ref{th: main: omega}.
This is the main reason for we set the orientation opposite to the usual one given in the above item, as a convention in this paper.
\end{rk}

\begin{rk}\label{rk: pre semi ribbon hyp}
Let us observe that if $\Gamma$ admits a pre-semi-ribbon graph $\Gamma_{\psr}$, it is necessarily a hyperbolic graph.
For, otherwise there exists a vertex $x\in \Gamma$ and a (local) sector $T$ centered at $x$ with the separatices with the same orientation.
By the definition of pre-semi-ribbon, there is no point of $\Gamma_{\psr}$ in $\Int(T)$.
But a smooth point of one of the separatices has an oriented sector intersecting $T$, which is contradiction to $\Gamma_{\psr}\cap\Int(T)=\emptyset$.
\end{rk}

In fact any hyperbolic graph admits a pre-semi-ribbon graph:

\begin{prop}\label{prop: pre semi ribbon}
Let $S$ be a Riemann surface with boundary and $\Gamma$ be an embedded graph.
Then $(S,\Gamma)$ admits a pre-semi-ribbon graph if and only if $\Gamma$ is a hyperbolic graph.
\end{prop}

The `only if' part was already observed in Remark \ref{rk: pre semi ribbon hyp}.
We postpone the proof of the `if' part to Appendix \ref{subsec: if part} because pre-semi-ribbon graph itself is \textit{not} important in the sequel and moreover the proof uses the idea of the proof of Theorem \ref{th: number of Eulerian paths}. 

\begin{defn}
A pre-semi-ribbon graph is called a \textit{semi-ribbon graph} if  $\Bd_{\sm}(\Gamma_{\psr})$ is connected.
A semi-ribbon graph of $\Gamma$ is denoted by $\Gamma_{\sr}$.
\end{defn}

Note that since $\Bd_{\sm}(\Gamma_{\sr})$ is connected, it is homeomorphic to the circle $S^{1}$.
This new constraint indeed provides a difference between pre-semi-ribbon and semi-ribbon graphs.
See Example \ref{non example} (2).

\begin{lem}\label{lem: conn hyp}
Let $S$ be a Riemann surface with boundary and $\Gamma$ be an embedded graph.
If $(S,\Gamma)$ admits a semi-ribbon graph $\Gamma_{\sr}$, $\Gamma$ is a path-connected hyperbolic graph.
\end{lem}

\proof
Since $\Bd_{\sm}(\Gamma_{\sr})\approx S^{1}$, it traverses every edge of $\Gamma$.
Hence $\Gamma$ is path-connected.
The hyperbolicity was already observed in Remark \ref{rk: pre semi ribbon hyp}.
\qed\vskip.3cm

\begin{example}\label{example: semi ribbon}
We give examples of semi-ribbon graphs.

(1)
A \textit{cactus} is a simply connected and connected finite union of
disks in $S^{2}$ such that any two disks intersect at most one point (Fig.\ \ref{Fig: 1-A}).
In other words its dual graph (by representing disks and intersection
points by vertices and edges respectively) becomes a tree.
The boundary of a cactus has the clockwise orientation induced from the one of
$S^{2}$ so that it becomes a finite oriented graph.
Since the boundary is deformed to an outer circle with clockwise orientation, it is equipped with a
semi-ribbon graph.

(2)
Let use consider a flow with slope $\frac23$ on the $2$-torus $T^{2}$
(Fig.\ \ref{Fig: 1-D}).
It has a semi-ribbon graph as in the figure after giving orientation from left-up to right-down.
\end{example}

\begin{example}\label{non example}
We already have a non-example in Fig.\ \ref{Fig: 1-C}.
Recall that this is non-hyperbolic so that it does not even admit pre-semi-ribbon graph.

Let us give an example of a pre-semi-graph which does not admit semi-graph.
Let us consider a hyperbolic graph in $S^{2}$ whose support (= graph without orientation) is given as in
Fig.\ \ref{Fig: 3-C}.
\begin{figure}[ht]
\begin{center}
\begin{tikzpicture}[line cap=round,line join=round,>=triangle 45,x=.7cm,y=.7cm]
\draw [rotate around={90.:(5.2,-1.2)}] (5.2,-1.2) ellipse (1.6 and 1.0583005244258357);
\draw [rotate around={0.:(5.2,-1.2)}] (5.2,-1.2) ellipse (3.939543120718437 and 1.6);
\end{tikzpicture}
\end{center}
\caption{}
\label{Fig: 3-C}
\end{figure}

This example does not admit a semi-ribbon graph but a pre-semi-ribbon graph for any hyperbolic orientation.
In any hyperbolic orientations, this has no semi-ribbon graph.
For, its complement in $S^{2}$ has $4$ connected components.
If a semi-ribbon graph exists, it lies in one of these $4$ components.
We leave it as an exercise that  these $4$ possibilities do not occur.

Recall that for any hyperbolic orientations, $\Gamma$ admits a pre-semi-ribbon graph due to the `if' part of Proposition \ref{prop: pre semi ribbon}.
\end{example}

We describe the semi-ribbon graphs in the genus $0$ Riemann surface.
This will correspond to the known description of $\omega(x)$ in genus $0$ (\cite[Theorem C]{JL}).
See Theorem \ref{th: JL}.

\begin{lem}\label{lem: geneus 0} 
Let $(S,\Gamma,\Gamma_{\sr})$ be a semi-ribbon graph for a connected
finite hyperbolic graph $\Gamma$ on $S=S^{2}$.
Then $(S,\Gamma,\Gamma_{\sr})$ is the semi-ribbon graph triple in Example \ref{example: semi ribbon} (1) (i.e., $\Gamma$ is a cactus boundary in $S^{2}$).
\end{lem}

\begin{figure}[ht]
\begin{center}
\begin{tikzpicture}[line cap=round,line join=round,>=triangle 45,x=.7cm,y=.7cm]
\draw (3.2,0.4)-- (5.6,0.4);
\draw (4.505,0.4) -- (4.4,0.265);
\draw (4.505,0.4) -- (4.4,0.535);
\draw (5.6,0.4)-- (8.,0.4);
\draw (6.905,0.4) -- (6.8,0.265);
\draw (6.905,0.4) -- (6.8,0.535);
\draw[->,decorate, decoration=zigzag] (5.6,-0.8)-- (5.6,-2.5);
\draw (3.2,-3.6)-- (8.,-3.6);
\draw (5.705,-3.6) -- (5.6,-3.735);
\draw (5.705,-3.6) -- (5.6,-3.465);
\begin{scriptsize}
\draw [fill=qqqqff] (3.2,0.4) circle (2.5pt);
\draw [fill=qqqqff] (5.6,0.4) circle (2.5pt);
\draw [color=black] (5.74,0.76) node {$v$};
\draw[color=black] (4.46,0.04) node {$e_1$};
\draw [fill=qqqqff] (8.,0.4) circle (2.5pt);
\draw[color=black] (6.86,0.04) node {$e_2$};
\draw [fill=qqqqff] (3.2,-3.6) circle (2.5pt);
\draw [fill=qqqqff] (8.,-3.6) circle (2.5pt);
\draw (4.2,-1.6) node {minimize};
\end{scriptsize}
\end{tikzpicture}
\end{center}
\caption{minimal graph}
\label{Fig: 3-e}
\end{figure}
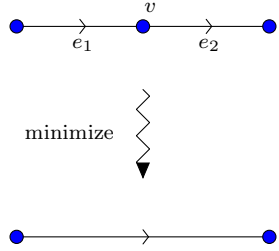

\proof
We replace $\Gamma$ into the minimal graph obtained in the following
steps:
Pick a vertex $v\in V$ having valency $2$ with  two distinct edges
$e_{1},e_{2}\in E$ satisfying $h(e_{1})=t(e_{2})=v$, if any.
Rule out $v$ from $V$ and $e_{1},e_{2}$ from $E$ and then add a new edge
amalgamating $e_{1},e_{2}$ (Fig.\ \ref{Fig: 3-e}).
By repeating this process for each vertex as above, the resulting
oriented graph is called the \textit{minimal} graph of $\Gamma$.
Note that due to the hyperbolicity assumption there is no vertex with
valency $2$ in the minimal graph unless $\Gamma$ is a loop with one vertex.
Note also that it suffices to check the lemma for the minimal graph.

If $\Gamma$ is a loop with one vertex, the lemma is obvious.
Indeed $\Gamma$ has clockwise or counter-clockwise orientation.
In the former case we are done by setting $\Bd_{\sm}(\Gamma_{\sr})$ to be an outer circle of $\Gamma$.
In the latter case we use an orientation reversing map of $S^{2}$.

Suppose that $\Gamma$ is not a loop.
This means that there is a vertex $v$ with valency $\ge4$.
Let us modify $\Gamma$ as follows:
Choose a sector not containing $\Bd_{\sm}(\Gamma)$ and then detach
$\Gamma$ along this sector.
We make the detached graph into the minimal graph and then denote it by
$\Gamma'$.
Note that $\Gamma'$ admits a natural pre-semi-ribbon graph induced from $\Gamma_{\sr}$. 
The above modification of $(\Gamma,\Gamma_{\sr})$ into $(\Gamma',\Gamma'_{\psr})$ is depicted as in Fig.\ \ref{Fig: 3-F}.
\begin{figure}[ht]
\begin{center}
\begin{tikzpicture}[line cap=round,line join=round,>=triangle 45,x=.8cm,y=.8cm]
\draw (2.,6.)-- (2.,4.4);
\draw (2.,5.095) -- (1.865,5.2);
\draw (2.,5.095) -- (2.135,5.2);
\draw (2.,4.4)-- (3.6,5.2);
\draw (2.8939148550549914,4.846957427527496) -- (2.8603738353924943,4.679252329215012);
\draw (2.8939148550549914,4.846957427527496) -- (2.739626164607506,4.92074767078499);
\draw (3.6,3.6)-- (2.,4.4);
\draw (2.7060851449450087,4.046957427527495) -- (2.860373835392494,4.120747670784988);
\draw (2.7060851449450087,4.046957427527495) -- (2.7396261646075053,3.879252329215011);
\draw (2.,4.4)-- (2.,2.8);
\draw (2.,3.495) -- (1.865,3.6);
\draw (2.,3.495) -- (2.135,3.6);
\draw (0.4,3.6)-- (2.,4.4);
\draw (1.293914855054991,4.046957427527495) -- (1.260373835392494,3.879252329215011);
\draw (1.293914855054991,4.046957427527495) -- (1.1396261646075054,4.120747670784988);
\draw (2.,4.4)-- (0.4,5.2);
\draw (1.1060851449450089,4.846957427527496) -- (1.2603738353924943,4.92074767078499);
\draw (1.1060851449450089,4.846957427527496) -- (1.1396261646075057,4.679252329215012);
\draw [->,decorate,decoration=zigzag] (4.4,4.4) -- (6.8,4.4);
\draw (8.8,6.)-- (8.8,4.4);
\draw (8.8,5.095) -- (8.665,5.2);
\draw (8.8,5.095) -- (8.935,5.2);
\draw (8.8,4.4)-- (8.8,2.8);
\draw (8.8,3.495) -- (8.665,3.6);
\draw (8.8,3.495) -- (8.935,3.6);
\draw (7.2,3.6)-- (8.8,4.4);
\draw (8.09391485505499,4.046957427527493) -- (8.060373835392493,3.879252329215009);
\draw (8.09391485505499,4.046957427527493) -- (7.939626164607504,4.120747670784986);
\draw (8.8,4.4)-- (7.2,5.2);
\draw (7.906085144945009,4.846957427527494) -- (8.060373835392495,4.920747670784987);
\draw (7.906085144945009,4.846957427527494) -- (7.939626164607506,4.67925232921501);
\draw (11.2,3.6)-- (10.,4.4);
\draw (10.512634719094525,4.058243520603649) -- (10.674884526490404,4.112326789735609);
\draw (10.512634719094525,4.058243520603649) -- (10.52511547350959,3.8876732102643907);
\draw (10.,4.4)-- (11.2,5.2);
\draw (10.687365280905473,4.858243520603649) -- (10.674884526490404,4.687673210264391);
\draw (10.687365280905473,4.858243520603649) -- (10.525115473509594,4.9123267897356095);
\draw [->,decorate,decoration=zigzag] (0.4,-0.4) -- (2.8,-0.4);
\draw (4.8,1.2)-- (4.8,-0.4);
\draw (4.8,0.295) -- (4.665,0.4);
\draw (4.8,0.295) -- (4.935,0.4);
\draw (4.8,-0.4)-- (4.8,-2.);
\draw (4.8,-1.305) -- (4.665,-1.2);
\draw (4.8,-1.305) -- (4.935,-1.2);
\draw (3.2,-1.2)-- (4.8,-0.4);
\draw (4.093914855054992,-0.7530425724725038) -- (4.060373835392495,-0.9207476707849883);
\draw (4.093914855054992,-0.7530425724725038) -- (3.9396261646075064,-0.6792523292150103);
\draw (4.8,-0.4)-- (3.2,0.4);
\draw (3.906085144945009,0.04695742752749541) -- (4.060373835392494,0.120747670784989);
\draw (3.906085144945009,0.04695742752749541) -- (3.939626164607506,-0.120747670784989);
\draw (7.2,-1.2)-- (6.,-0.4);
\draw (6.512634719094527,-0.7417564793963497) -- (6.674884526490407,-0.6876732102643898);
\draw (6.512634719094527,-0.7417564793963497) -- (6.5251154735095955,-0.9123267897356078);
\draw (6.,-0.4)-- (7.2,0.4);
\draw (6.6873652809054756,0.058243520603649616) -- (6.674884526490407,-0.11232678973560843);
\draw (6.6873652809054756,0.058243520603649616) -- (6.5251154735095955,0.11232678973560957);
\draw [->,decorate,decoration=zigzag] (8.,-0.4) -- (10.4,-0.4);
\draw (12.4,1.2)-- (12.4,-0.4);
\draw (12.4,0.295) -- (12.265,0.4);
\draw (12.4,0.295) -- (12.535,0.4);
\draw (12.4,-0.4)-- (12.4,-2.);
\draw (12.4,-1.305) -- (12.265,-1.2);
\draw (12.4,-1.305) -- (12.535,-1.2);
\draw (10.8,-1.2)-- (12.4,-0.4);
\draw (11.69391485505499,-0.7530425724725038) -- (11.660373835392493,-0.9207476707849883);
\draw (11.69391485505499,-0.7530425724725038) -- (11.539626164607503,-0.6792523292150103);
\draw (12.4,-0.4)-- (10.8,0.4);
\draw (11.50608514494501,0.04695742752749541) -- (11.660373835392493,0.120747670784989);
\draw (11.50608514494501,0.04695742752749541) -- (11.539626164607503,-0.120747670784989);
\draw (14.8,-1.2)-- (13.6,-0.4);
\draw (14.112634719094526,-0.7417564793963508) -- (14.274884526490405,-0.6876732102643909);
\draw (14.112634719094526,-0.7417564793963508) -- (14.125115473509592,-0.9123267897356089);
\draw (13.6,-0.4)-- (14.8,0.4);
\draw (14.287365280905474,0.058243520603649616) -- (14.274884526490405,-0.11232678973560843);
\draw (14.287365280905474,0.058243520603649616) -- (14.125115473509595,0.11232678973560957);
\draw[color=gray,arrow data={.3}{stealth},arrow data={.7}{stealth}] plot [color=gray,smooth, tension=.7] coordinates
{(2.4,6.) (2.4,5.36) (2.8,5.2) (3.6,5.6)};
\draw[color=gray,arrow data={.4}{stealth},arrow data={.7}{stealth}] plot [color=gray,smooth, tension=.7] coordinates
{(0.4,4.) (1.02,4.4) (0.4,4.8)};
\draw[color=gray,arrow data={.4}{stealth},arrow data={.7}{stealth}] plot [color=gray,smooth, tension=.7] coordinates
{(3.6,3.2) (2.8,3.6) (2.38,3.4) (2.4,2.8)};
\draw[color=gray,arrow data={.3}{stealth},arrow data={.7}{stealth}] plot [color=gray,smooth, tension=.2] coordinates
{(12.8,1.2) (12.8,-2.)};
\draw[color=gray,arrow data={.4}{stealth},arrow data={.7}{stealth}] plot [color=gray,smooth, tension=.7] coordinates
{(10.8,-0.8) (11.4,-0.4) (10.8,0.)};
\draw[color=gray,arrow data={.3}{stealth},arrow data={.7}{stealth}] plot [color=gray,smooth, tension=.7] coordinates
{(14.8,-1.6) (13.6,-0.8) (13.3,-0.4) (13.6,0.) (14.8,0.8)};
\begin{scriptsize}
\draw [fill=qqqqff] (2.,6.) circle (2.5pt);
\draw [fill=qqqqff] (2.,4.4) circle (2.5pt);
\draw [fill=qqqqff] (3.6,5.2) circle (2.5pt);
\draw [fill=qqqqff] (3.6,3.6) circle (2.5pt);
\draw [fill=qqqqff] (2.,2.8) circle (2.5pt);
\draw [fill=qqqqff] (0.4,3.6) circle (2.5pt);
\draw [fill=qqqqff] (0.4,5.2) circle (2.5pt);
\draw[color=black] (5.46,4.94) node {detach};
\draw [fill=qqqqff] (8.8,6.) circle (2.5pt);
\draw [fill=qqqqff] (8.8,4.4) circle (2.5pt);
\draw [fill=qqqqff] (11.2,5.2) circle (2.5pt);
\draw [fill=qqqqff] (11.2,3.6) circle (2.5pt);
\draw [fill=qqqqff] (8.8,2.8) circle (2.5pt);
\draw [fill=qqqqff] (7.2,3.6) circle (2.5pt);
\draw [fill=qqqqff] (7.2,5.2) circle (2.5pt);
\draw [fill=qqqqff] (10.,4.4) circle (2.5pt);
\draw[color=black] (1.5,0.2) node {minimize};
\draw [fill=qqqqff] (4.8,1.2) circle (2.5pt);
\draw [fill=qqqqff] (4.8,-0.4) circle (2.5pt);
\draw [fill=qqqqff] (7.2,0.4) circle (2.5pt);
\draw [fill=qqqqff] (7.2,-1.2) circle (2.5pt);
\draw [fill=qqqqff] (4.8,-2.) circle (2.5pt);
\draw [fill=qqqqff] (3.2,-1.2) circle (2.5pt);
\draw [fill=qqqqff] (3.2,0.4) circle (2.5pt);
\draw[color=black] (9.08,0.2) node {set $\Bd_{\sm}(\Gamma'_{\psr})$};
\draw [fill=qqqqff] (12.4,1.2) circle (2.5pt);
\draw [fill=qqqqff] (12.4,-0.4) circle (2.5pt);
\draw [fill=qqqqff] (14.8,0.4) circle (2.5pt);
\draw [fill=qqqqff] (14.8,-1.2) circle (2.5pt);
\draw [fill=qqqqff] (12.4,-2.) circle (2.5pt);
\draw [fill=qqqqff] (10.8,-1.2) circle (2.5pt);
\draw [fill=qqqqff] (10.8,0.4) circle (2.5pt); 
\end{scriptsize}
\end{tikzpicture}
\end{center}
\caption{}
\label{Fig: 3-F}
\end{figure}

Note that $\Gamma'$ has at most two connected components.
I.e., it depends on if there is a (not necessarily oriented) path between the detached vertices.
Each component is equipped with a semi-ribbon graph by the construction.

We claim $\Gamma'$ is disconnected. 
Recall that $\Bd_{\sm}(\Gamma_{\sr})$ changes into $\Bd_{\sm}(\Gamma'_{\psr})$ by replacing the left crossing to the right one in Fig.\ \ref{Fig: 3-g}.
\begin{figure}[ht]
\begin{center}
\begin{tikzpicture}[line cap=round,line join=round,>=triangle 45,x=1.1cm,y=1.1cm]
\draw [color=gray,arrow data={.5}{stealth},shift={(-0.4,2.)}] plot[domain=4.124386376837123:5.3003915839322575,variable=\t]({1.*1.4422205101855954*cos(\t r)+0.*1.4422205101855954*sin(\t r)},{0.*1.4422205101855954*cos(\t r)+1.*1.4422205101855954*sin(\t r)});
\draw [color=gray,arrow data={.5}{stealth},shift={(-0.4,-2.)}] plot[domain=0.9827937232473286:2.158798930342464,variable=\t]({1.*1.4422205101855963*cos(\t r)+0.*1.4422205101855963*sin(\t r)},{0.*1.4422205101855963*cos(\t r)+1.*1.4422205101855963*sin(\t r)});
\draw [->,decorate,decoration=zigzag] (1.6,0.) -- (3.6,0.);
\draw [color=gray,arrow data={.5}{stealth},shift={(3.6,0.)}] plot[domain=0.5880026035475675:-0.5880026035475678,variable=\t]({1.*1.4422205101855956*cos(\t r)+0.*1.4422205101855956*sin(\t r)},{0.*1.4422205101855956*cos(\t r)+1.*1.4422205101855956*sin(\t r)});
\draw [color=gray,arrow data={.5}{stealth},shift={(7.6,0.)}] plot[domain=3.7295952571373605:2.5535900500422257,variable=\t]({1.*1.442220510185596*cos(\t r)+0.*1.442220510185596*sin(\t r)},{0.*1.442220510185596*cos(\t r)+1.*1.442220510185596*sin(\t r)});
\begin{scriptsize}
\draw[color=black] (-1.46,-0.84) node {$d$};
\draw[color=black,shift={(6,0)}] (-1.46,-0.84) node {$d$};
\draw[color=black] (-1.46,0.96) node {$a$};
\draw[color=black,shift={(6,0)}] (-1.46,0.96) node {$a$};
\draw[color=black] (0.64,0.96) node {$b$};
\draw[color=black,shift={(6,0)}] (0.64,0.96) node {$b$};
\draw[color=black] (0.64,-0.84) node {$c$};
\draw[color=black,shift={(6,0)}] (0.64,-0.84) node {$c$};
\end{scriptsize}
\end{tikzpicture}
\end{center}
\caption{}
\label{Fig: 3-g}
\end{figure}

We may assume the orientation of $\Bd_{\sm}(\Gamma)\approx S^{1}$ traverses $S^{1}$ in the cyclic order
	$$
	a\to b\to c\to d\to a.
	$$
On the other hand the orientation  of $\Bd_{\sm}(\Gamma'_{\psr})$ traverses  in the cyclic orders
	\begin{equation}\label{eq: cyclic orders}
	a\to d\to a,\quad c\to b\to c.
	\end{equation}
Therefore $\Bd_{\sm}(\Gamma'_{\psr})$ consists of two copies of $S^{1}$ and thus $\Gamma'$ is disconnected.

We denote the two connected components of $\Gamma'$ by $\Gamma'_{1},\Gamma'_{2}$.
As a consequence $\Gamma'_{\psr}$ has two connected components and each component becomes a semi-ribbon graph $(\Gamma_{i}')_{\sr}$ of $\Gamma'_{i}$ ($i=1,2$).
By induction on $\#E$, both $\Gamma'_{1},\Gamma'_{2}$ are cactus boundaries with the same orientation due to the cyclic orders \eqref{eq: cyclic orders}.
And $\Bd_{\sm}((\Gamma'_{i})_{\sr})$ encloses the cactus $(i=1,2)$. 
Looking at the fourth graph of Fig.\ \ref{Fig: 3-F}, the original graph $\Gamma$ is the boundary of the cactus obtained by pasting the two cacti along the detached vertices. 
This identifies $(S,\Gamma,\Gamma_{\sr})$ as in the lemma.
\qed\vskip.3cm


\subsection{Eulerian path}\label{subsec: Euler}

A path in a connected oriented graph $\Gamma$ is called \textit{Eulerian} if each edge in the path is used exactly once.
A \textit{full} Eulerian path of $\Gamma$ is an Eulerian path using every edge of $\Gamma$ exactly once.
$\Gamma$ itself is called \textit{Eulerian} if it admits a full Eulerian path.
$\Gamma$ is Eulerian if and only if $\# t\inv(v)=\# h\inv(v)$ for each vertex $v$.
The number of full Eulerian paths with a fixed starting point is explicitly
given for a connected finite oriented graph (a.k.a.\ the BEST theorem in
graph theory).

On the other hand if a connected finite oriented $\Gamma$ is equipped
with a semi-ribbon graph $(S,\Gamma,\Gamma_{\sr})$, the orientation of
$\Bd_{\sm}(\Gamma_{\sr})$ gives one full Eulerian path on $\Gamma$.
However \textit{not} every full Eulerian path on a given connected finite
hyperbolic graph $\Gamma$ comes from choice of $(S,\Gamma_{\sr})$.
We do not pursue this general question in this paper, but give a partial
answer only when $S$ is orientable.

\begin{thm}\label{th: number of Eulerian paths}
Let $(S,\Gamma,\Gamma_\sr)$ be a semi-ribbon graph, where $S$ is a connected Riemann surface with boundary and $\Gamma$ is a connected finite hyperbolic graph.
Then we have the following description of $(S,\Gamma,\Gamma_\sr)$.
	\begin{itemize}
	\item[(i)] $S$ is  a connected sum of $S^{2}$ and a Riemann surface with boundary along finitely many punctures disjoint from $\Gamma$.
	\item[(ii)] There is a cactus $D$ in $S^{2}$ such that $\Gamma=\Bd(D)$ with clockwise orientation and $\Bd_{\sm}(\Gamma_{\sr})\, (\approx S^{1})$ encloses $D$.
	\end{itemize}

In particular the number of full Eulerian paths from semi-ribbon graphs on Riemann surfaces is $1$.
\end{thm}

The proof will appear in the next subsection.

\begin{rk}
Let us compare the statements in Theorems \ref{th: main: analytic version} and \ref{th: number of Eulerian paths}.

(1)
In the latter theorem the Riemann surfaces may have boundary while in the former, without boundary.
As explained in \S\ref{subsec: main theorem}, the former theorem will be generalized to the Riemann surfaces with boundary in Theorem \ref{th: main: omega}.

(2)
In the former theorem there is a statement on convergence of flow.
This part is replaced by semi-ribbon graph in the latter theorem.
In other words the enclosed area between $\Bd_{\sm}(\Gamma_{\sr})$ and $\Bd_{\sing}(\Gamma_{\sr})\, (=\Gamma)$ is attracted to $\Gamma$.
This phenomenon has been a key fact in the topological study of flows on genus 0 surfaces in various literatures. 

In the latter theorem there is also additional information on orientation.
Since $\omega(x)$ consists of finite orbits (but not a single point), it has the orientation from flow direction.
This corresponds to the orientation of $\Gamma$.

(3)
There is a further notice on orientation in the latter theorem.
$\Gamma$ has only clockwise orientation on $S^{2}$ which is due to the orientation sector rule in the definition of pre-semi-ribbon graph.
However the counter-clockwise orientation of $\Gamma$ is also included in the latter theorem by reversing the orientation used in the connected sum.
\end{rk}

A \textit{simple} path of $\Gamma$ means a path without self-intersection.
For any given path one can construct a simple path with the same starting and terminal points:
For a self-intersection vertex on the path, we delete the loops attached to it.
We call this procedure \textit{simplification}. 

\begin{lem}\label{lem: existence of minimal closed path}
Let $S$ be a Riemann surface with boundary.
Let $\Gamma$ be a connected finite graph in $S$ such that $\#h\inv(v),\ \#t\inv(v)\ge1$ for each vertex $v\in V$.
Let $D$ be an embedded disk in $S$ such that  $\Bd(D)$ is a path in $\Gamma$ where the orientation of $\Bd(D)$ is induced from $D$.
If $\Int(D)\cap \Gamma\neq\emptyset$, there exists a simple closed path of $\Gamma$ in $D$ not equal to $\Bd(D)$.
\end{lem}

\proof
By the assumption there is at least one edge in $D$ not contained in $\Bd(D)$, say $e$.
Suppose there is no vertex of $\Gamma$ in $\Int(D)$.
Then $e$ is a path from $p$ to $q$ for some $p,q\in \Bd(D)$.
If $p=q$, then $e$ is a loop and thus this loop satisfies the lemma.
If $p\neq q$, $e$ separates $D$ into two connected components, say $D_{1},D_{2}$.
One of $\Bd(D_{1}),\Bd(D_{2})$ with the orientations from $D_{1},D_{2}$ is a simple closed path of $\Gamma$.
It is clear that this path satisfies the lemma.

Suppose next that there is a vertex in $\Int(D)$.
Then there is an edge $e$ satisfying \textit{either}
	\begin{itemize}
	\item[(i)]
	$t(e)\in \Bd(D),\ h(e)\in\Int(D)$ \textit{or}
	\item[(ii)]
	$h(e)\in \Bd(D),\ t(e)\in\Int(D)$
	\end{itemize}
(see Fig.\ \ref{Fig: 3-H}).
\begin{figure}[ht]
\begin{center}
\begin{tikzpicture}[line cap=round,line join=round,>=triangle 45,x=.7cm,y=.7cm]
\draw(0.,0.4) circle (2.8);
\draw(8.8,0.4) circle (2.8);
\draw (0.,3.2)-- (0.,0.8);
\draw (0.,1.895) -- (-0.135,2.);
\draw (0.,1.895) -- (0.135,2.);
\draw (8.8,0.8)-- (8.8,3.2);
\draw (8.8,2.105) -- (8.935,2.);
\draw (8.8,2.105) -- (8.665,2.);
\begin{scriptsize}
\draw[color=black] (-2.7,4.02) node {(i)};
\draw[color=black] (6.14,3.96) node {(ii)};
\draw [fill=qqqqff] (0.,3.2) circle (2.5pt);
\draw[color=black] (0.14,3.56) node {$t(e)$};
\draw [fill=qqqqff] (0.,0.8) circle (2.5pt);
\draw[color=black] (0.14,0.44) node {$h(e)$};
\draw[color=black] (-0.26,2.16) node {$e$};
\draw [fill=qqqqff] (8.8,0.8) circle (2.5pt);
\draw[color=black] (8.94,0.44) node {$t(e)$};
\draw [fill=qqqqff] (8.8,3.2) circle (2.5pt);
\draw[color=black] (8.94,3.56) node {$h(e)$};
\draw[color=black] (9.18,2.16) node {$e$};
\end{scriptsize}
\end{tikzpicture}
\end{center}
\caption{}
\label{Fig: 3-H}
\end{figure}

Let us assume (i) first.
We consider a simple path $P$ from $v:=h(e)$ in $D$ with the maximal number of edges.
Note that by the assumption $\#t\inv(v)\ge1$, there is at least one path from $v$.
Let $v'$ be the terminal vertex of $P$.
Suppose first that $v'\in \Bd(D)$.
Let $C'$ be the subpath in $\Bd(D)$ from $v'$ to $t(e)$.
Now the closed path formed by $C',P$ and $e$ satisfies the lemma after simplification if necessary.

Even in the case when $P\cap\Bd(D)\neq \emptyset$, a similar construction works by replacing $v'$ into any intersection vertex in $P\cap \Bd(D)$.
We omit the details.

The remaining case is $P\subset \Int(D)$.
By the maximality of number of edges and the assumption $\#t\inv(v')\ge1$, the longer path $P'$ obtained by adding to $P$ and an edge $e'$ with $t(e')=v'$ has a self-intersection vertex, say $v''$.
The subpath of $P'$ from $v''$ to itself satisfies the lemma. 

We complete the proof of the lemma in the case (i).
The case (ii) is similar.
\qed\vskip.3cm

Note that the condition $\#h\inv(v),\, \#t\inv(v)\ge1$ of the above lemma is weaker than hyperbolicity.
In fact any finite graph with this condition allows a simple closed path. 
This is purely graph-theoretic and irrelevant to the ambient surfaces.
See Lemma \ref{lem: existence of simple closed path}.

Now we impose hyperbolicity of $\Gamma$.

\begin{lem}\label{lem: homology 0}
Let $\Gamma$ be a connected finite hyperbolic graph in a Riemann surface $S$ with boundary.
If any simple closed path has class $0$ in $H_{1}(S,\zz)$, there exists an embedded disk $D$ in $S$ containing $\Gamma$.
\end{lem}

\proof
We take any simple closed path of $\Gamma$ which is obtained by simplification of a full Eulerian path of $\Gamma$.
By the assumption on the trivial homology class, the path is $\Bd(D_{0})$ for some embedded disk $D_{0}$ in $S$.
By Lemma \ref{lem: existence of minimal closed path}, we may assume that $\Int(D_{0})$ contains no edge of $\Gamma$.
So any edge of $\Gamma$ not contained in $\Bd(D_{0})$ lies outside $D_{0}$.
The local picture of $\Gamma$ near $D_{0}$ looks as in left graph of Fig.\ \ref{Fig: 3-I}.
\begin{figure}[ht]
\begin{center}
\begin{tikzpicture}[line cap=round,line join=round,>=triangle 45,x=.7cm,y=.7cm]
\draw (-1.6,3.6)-- (1.6,3.6);
\draw (0.105,3.6) -- (0.,3.465);
\draw (0.105,3.6) -- (0.,3.735);
\draw (1.6,3.6)-- (1.6,0.4);
\draw (1.6,1.895) -- (1.465,2.);
\draw (1.6,1.895) -- (1.735,2.);
\draw (1.6,0.4)-- (-1.6,0.4);
\draw (-0.105,0.4) -- (0.,0.535);
\draw (-0.105,0.4) -- (0.,0.265);
\draw (-1.6,0.4)-- (-1.6,3.6);
\draw (-1.6,2.105) -- (-1.465,2.);
\draw (-1.6,2.105) -- (-1.735,2.);
\draw (1.6,3.6)-- (1.6,5.2);
\draw (1.6,4.505) -- (1.735,4.4);
\draw (1.6,4.505) -- (1.465,4.4);
\draw (3.2,3.6)-- (1.6,3.6);
\draw (2.295,3.6) -- (2.4,3.735);
\draw (2.295,3.6) -- (2.4,3.465);
\draw (-3.2,0.4)-- (-1.6,0.4);
\draw (-2.295,0.4) -- (-2.4,0.265);
\draw (-2.295,0.4) -- (-2.4,0.535);
\draw (-1.6,0.4)-- (-1.6,-1.2);
\draw (-1.6,-0.505) -- (-1.735,-0.4);
\draw (-1.6,-0.505) -- (-1.465,-0.4);
\draw [->,decorate,decoration=zigzag] (4.,2.) -- (6.4,2.);
\draw (9.6,2.)-- (9.6,3.6);
\draw (9.6,2.905) -- (9.735,2.8);
\draw (9.6,2.905) -- (9.465,2.8);
\draw (11.2,2.)-- (9.6,2.);
\draw (10.295,2.) -- (10.4,2.135);
\draw (10.295,2.) -- (10.4,1.865);
\draw (9.6,2.)-- (9.6,0.);
\draw (9.6,0.895) -- (9.465,1.);
\draw (9.6,0.895) -- (9.735,1.);
\draw (8.,2.)-- (9.6,2.);
\draw (8.905,2.) -- (8.8,1.865);
\draw (8.905,2.) -- (8.8,2.135);
\begin{scriptsize}
\draw [fill=qqqqff] (-1.6,3.6) circle (2.5pt);
\draw [fill=qqqqff] (1.6,3.6) circle (2.5pt);
\draw [fill=qqqqff] (1.6,0.4) circle (2.5pt);
\draw [fill=qqqqff] (-1.6,0.4) circle (2.5pt);
\draw[color=black] (2.08,4.56) node {$e_1$};
\draw[color=black] (2.46,4.00) node {$e_2$};
\draw[color=black] (-2.34,0.04) node {$e_3$};
\draw[color=black] (-2.06,-0.44) node {$e_4$};
\draw[color=black] (5.,2.6) node {shrink $D_{0}$};
\draw [fill=qqqqff] (9.6,2.) circle (2.5pt);
\draw[color=black] (10.08,2.96) node {$e_1$};
\draw[color=black] (10.46,2.40) node {$e_2$};
\draw[color=black] (9.14,0.96) node {$e_3$};
\draw[color=black] (8.86,1.64) node {$e_4$};
\draw[color=black] (0.,2.04) node {$D_{0}$};
\end{scriptsize}
\end{tikzpicture}
\end{center}
\caption{}
\label{Fig: 3-I}
\end{figure}

If $\Gamma=\Bd(D_{0})$, $D_{0}$ itself satisfies the lemma. 
Thus we assume there are edges outside $D_{0}$.
We shrink $D_{0}$ to a point.
This process induces $(\oS,\oGamma)$ from $(S,\Gamma)$.
To be precise $\oS$ is obtained by $S$ by shrinking the disk $D_{0}$.
Thus $\oS$ is homeomorphic to $S$.
$\oGamma$ is obtained by replacing the edges and the vertices in $\Bd(D_{0})$ into one vertex.
This shrunken vertex is the central vertex in the right graph in Fig.\ \ref{Fig: 3-I}.
It is a connected hyperbolic graph with fewer edges as one can see the right graph in Fig.\ \ref{Fig: 3-I}.
It is clear that any simple closed path of $\oGamma$ has also class $0$ in $H_{1}(\oS,\zz)$.
Thus $(\oS,\oGamma)$ satisfies the induction hypothesis on $\#E$. 
Therefore there is an embedded disk $\oD$ containing $\oGamma$.

Now the reverse process of shrinking $D_{0}$ expands $\oD$ to an embedded disk in $S$.
This disk contains $\Gamma$, which proves the lemma. 
\qed\vskip.3cm

The following assertion will not be used in the proof of Theorem \ref{th: number of Eulerian paths}.
It will be used only in the appendix. 
So the readers can skip, but we record here the hyperbolicity of pull-back graph under cutting surgery.

\begin{prop}\label{prop: }
Let $\Gamma$ be a connected finite hyperbolic graph in a Riemann surface $S$ with boundary.
Let $C$ be a simple closed path of $\Gamma$.
For $\pi\colon \tS\to S$ the cutting surgery along $C$, $\pi\inv(\Gamma)$ is a finite hyperbolic graph.
\end{prop}

\proof
This comes from a local picture at a vertex in $C$.
See Fig.\ \ref{Fig: 3-K}.
\begin{figure}[ht]
\begin{center}
\begin{tikzpicture}[line cap=round,line join=round,>=triangle 45,x=.7cm,y=.7cm]
\draw (-0.4,-0.8)-- (-0.4,1.6);
\draw (-0.4,0.505) -- (-0.265,0.4);
\draw (-0.4,0.505) -- (-0.535,0.4);
\draw (-0.4,1.6)-- (-0.4,4.);
\draw (-0.4,2.905) -- (-0.265,2.8);
\draw (-0.4,2.905) -- (-0.535,2.8);
\draw (0.4,-0.8)-- (0.4,1.6);
\draw (0.4,0.505) -- (0.535,0.4);
\draw (0.4,0.505) -- (0.265,0.4);
\draw (0.4,1.6)-- (0.4,4.);
\draw (0.4,2.905) -- (0.535,2.8);
\draw (0.4,2.905) -- (0.265,2.8);
\draw (-0.4,1.6)-- (-1.6,-0.4);
\draw (-1.0540220543198904,0.5099632428001828) -- (-1.1157615449711937,0.6694569269827164);
\draw (-1.0540220543198904,0.5099632428001828) -- (-0.8842384550288068,0.5305430730172841);
\draw (-2.4,0.8)-- (-0.4,1.6);
\draw (-1.3025099474570476,1.238996021017181) -- (-1.349862258692196,1.0746556467304906);
\draw (-1.3025099474570476,1.238996021017181) -- (-1.450137741307804,1.3253443532695104);
\draw (-0.4,1.6)-- (-2.4,2.4);
\draw (-1.4974900525429524,2.0389960210171805) -- (-1.349862258692196,2.1253443532695093);
\draw (-1.4974900525429524,2.0389960210171805) -- (-1.450137741307804,1.8746556467304891);
\draw (-1.6,3.6)-- (-0.4,1.6);
\draw (-0.9459779456801095,2.5099632428001826) -- (-1.1157615449711933,2.5305430730172835);
\draw (-0.9459779456801095,2.5099632428001826) -- (-0.8842384550288063,2.6694569269827153);
\draw (0.4,1.6)-- (2.,0.4);
\draw (1.284,0.937) -- (1.119,0.892);
\draw (1.284,0.937) -- (1.281,1.108);
\draw (2.,2.8)-- (0.4,1.6);
\draw (1.116,2.137) -- (1.119,2.308);
\draw (1.116,2.137) -- (1.281,2.092);
\draw [->] (3.2,1.6) -- (5.6,1.6);
\draw (8.8,-0.8)-- (8.8,1.6);
\draw (8.8,0.505) -- (8.935,0.4);
\draw (8.8,0.505) -- (8.665,0.4);
\draw (8.8,1.6)-- (8.8,4.);
\draw (8.8,2.905) -- (8.935,2.8);
\draw (8.8,2.905) -- (8.665,2.8);
\draw (8.8,1.6)-- (7.6,-0.4);
\draw (8.14597794568011,0.5099632428001839) -- (8.084238455028808,0.6694569269827175);
\draw (8.14597794568011,0.5099632428001839) -- (8.315761544971194,0.5305430730172852);
\draw (6.8,0.8)-- (8.8,1.6);
\draw (7.897490052542953,1.2389960210171824) -- (7.850137741307804,1.0746556467304917);
\draw (7.897490052542953,1.2389960210171824) -- (7.749862258692197,1.3253443532695115);
\draw (8.8,1.6)-- (6.8,2.4);
\draw (7.702509947457047,2.0389960210171827) -- (7.850137741307804,2.1253443532695115);
\draw (7.702509947457047,2.0389960210171827) -- (7.749862258692194,1.8746556467304913);
\draw (7.6,3.6)-- (8.8,1.6);
\draw (8.25402205431989,2.509963242800184) -- (8.084238455028807,2.5305430730172853);
\draw (8.25402205431989,2.509963242800184) -- (8.315761544971192,2.669456926982717);
\draw (10.4,2.8)-- (8.8,1.6);
\draw (9.516,2.137) -- (9.519,2.308);
\draw (9.516,2.137) -- (9.681,2.092);
\draw (8.8,1.6)-- (10.4,0.4);
\draw (9.684,0.937) -- (9.519,0.892);
\draw (9.684,0.937) -- (9.681,1.108);
\begin{scriptsize}
\draw [fill=qqqqff] (-0.4,1.6) circle (2.5pt);
\draw [fill=qqqqff] (0.4,1.6) circle (2.5pt);
\draw [fill=qqqqff] (8.8,1.6) circle (2.5pt);
\draw[color=black] (4.46,1.84) node {$\pi$};
\end{scriptsize}
\end{tikzpicture}
\end{center}
\caption{}
\label{Fig: 3-K}
\end{figure}

We explain the details in Fig.\ \ref{Fig: 3-K}.
In the left graph, there are two groups of vertical edges.
The group in the left (resp.\ right) lies in $C_1$ (resp.\ $C_2$) in $\pi\inv(C)$.
The left graph maps to the right one via $\pi$ by identifying upper (resp.\ lower) vertical arrows.
This figure shows that hyperbolicity property at the central vertex is invariant under the cutting surgery.
\qed\vskip.3cm

We finish this section by giving the proof of the main theorem.
\vskip.3cm


\subsection{Proof of Theorem \ref{th: number of Eulerian paths}} 
This subsection is devoted to the proof of Theorem \ref{th: number of Eulerian paths}.

Suppose any simple closed path of $\Gamma$ has the trivial class in $H_{1}(S,\zz)$.
By Lemma \ref{lem: homology 0} there is an embedded disk $D$ in $S$ containing $\Gamma$.
After replacing into a slightly bigger disk if necessary, we may assume $\Gamma_{\sr}\subset D$.
By capping-off $\Bd(D)$ by a disk, $\Gamma_{\sr}$ is a semi-ribbon graph embedded in $S^{2}$. 
Thus by the genus 0 (Lemma \ref{lem: geneus 0}), $(S^{2},\Gamma,\Gamma_{\sr})$ is given by a cactus on $S^{2}$.
There are the two cases due to the pre-semi-ribbon graph axiom: 
	\begin{enumerate}
	\item
 	 $\Bd_{\sm}(\Gamma_{\sr})$ encloses $\Gamma$ in $D$,  
	\item $\Gamma$ is a circle enclosing $\Bd_{\sm}(\Gamma_{\sr})$ in $D$.
	\end{enumerate}
In each case we will complete the proof of the theorem.
In the case (1), by replacing the above capping-off disk to $S\setminus D$, $S$ is a connected sum of $S^{2}$ and $S$ itself along the puncture $\Bd(D)$.
Similarly in the case (2), $S$ is a connected sum of $S^{2}$ and $S$ itself along the puncture $\Bd(D)$.

We assume now that there is a simple closed path $C$ of $\Gamma$ with $[C]\neq0$ in $H_{1}(S,\zz)$.
Let $\pi\colon \tS\to S$ be the cutting surgery along $C$.
Since $\pi\inv(\Bd_{\sm}(\Gamma_{\sr}))$ has the homeomorphic image $\Bd_{\sm}(\Gamma_{\sr})\, (\approx S^{1})$ via $\pi$, it defines an Eulerian path of $\pi\inv(\Gamma)$.
We take the subgraph $\tGamma$ of $\pi\inv(\Gamma)$ which has the above Eulerian path as a full Eulerian path.
Now there is the induced semi-ribbon-graph $\tGamma_{\sr}$ with the smooth boundary $\pi\inv(\Bd_{\sm}(\Gamma_{\sr}))$. 
To see this, let $U_{x}$ be an open neighborhood of $x\in \Gamma$ in $S$ satisfying the axiom of the pre-semi-ribbon graph.
Let $\tx\in \tGamma$ be an inverse image of $x$.
Let $U_{\tx}$ be the connected component of $\pi\inv(U_{x})$ which intersects $\pi\inv(\Bd_{\sm}(\Gamma_{\sr}))$.
These open subsets of $\tS$ satisfy the axiom of the pre-semi-ribbon graph.
Thus we have the induced semi-ribbon graph $\tGamma_{\sr}$.
By Lemma \ref{lem: conn hyp}, $\tGamma$ is a connected hyperbolic graph.

We proceed the proof in the two cases: $\tS$ is connected or not.
Suppose first that $\tS$ is disconnected.
It has two connected components $\tS_{1},\tS_{2}$ with $g(\tS_{1}),g(\tS_{2})<g(S)$ by Lemma \ref{lem: genus: disconnected case}.
Since $\tGamma_{\sr}$ is connected, it lies in one of $\tS_{1},\tS_{2}$, say $\tS_{1}$.
By genus induction, $(\tS_{1},\tGamma,\tGamma_{\sr})$ satisfies the theorem.
Hence $(S,\Gamma,\Gamma_{\sr})$ is constructed from $(\tS_{1},\tGamma,\tGamma_{\sr})$ by the usual connected sum of $\tS_{1},\tS_{2}$ along $C_{1},C_{2}$.
Thus we have the description of the theorem in this case.

If $\tS$ is connected, by Lemma \ref{lem: genus: connected case}, $g(\tS)=g(S)-1$.
Let $\tC_{i}:=C_{i}\cap \tGamma$, $i=1,2$, subgraphs of $\tGamma$.
We claim that \textit{either} $\tC_{1}=C_{1},\ \tC_{2}=\emptyset$
\textit{or} $\tC_{1}=\emptyset,\ \tC_{2}=C_{2}$.
We prove this claim in two steps: 
First either $\tC_{1}$ or $\tC_{2}$ contains no edge of $\tGamma$. 
Secondly if there is no edge of $\tGamma$ in $\tC_{2}$, then no vertex in $\tC_{2}$ i.e.\ $\tC_{2}=\emptyset$.

We check the first step.
Suppose the contrary: both $\tC_{1},\tC_{2}$ have edges.
Note that by the construction of $\tGamma$, the edges of $\tC_{1},\tC_{2}$ map bijectively to those of $C$ via $\pi$.
Therefore there are edges $\widetilde{e}_1,\widetilde{e}_2$ in $\tC_1,\tC_2$ respectively such that their $\pi$-images $e_1:=\pi(\widetilde{e}_1)$ and $e_2:=\pi(\widetilde{e}_2)$ have a common vertex $v$.
Locally at $\pi\inv(v)=\{\tv_{1},\tv_{2}\}$,  $\tGamma$ is depicted as in Fig.\ \ref{Fig: 3-L}.
\begin{figure}[ht]
\begin{center}
\begin{tikzpicture}[line cap=round,line join=round,>=triangle 45,x=.7cm,y=.7cm]
\draw (-0.4,-0.4)-- (-0.4,1.6);
\draw (-0.4,0.705) -- (-0.265,0.6);
\draw (-0.4,0.705) -- (-0.535,0.6);
\draw (-0.4,1.6)-- (-0.4,3.6);
\draw (-0.4,2.705) -- (-0.265,2.6);
\draw (-0.4,2.705) -- (-0.535,2.6);
\draw (-1.6,3.2)-- (-0.4,1.6);
\draw (-0.937,2.316) -- (-1.108,2.319);
\draw (-0.937,2.316) -- (-0.892,2.481);
\draw (-0.4,1.6)-- (-2.4,2.4);
\draw (-1.4974900525429524,2.038996021017181) -- (-1.349862258692196,2.12534435326951);
\draw (-1.4974900525429524,2.038996021017181) -- (-1.450137741307804,1.8746556467304898);
\draw (-2.4,0.8)-- (-0.4,1.6);
\draw (-1.3025099474570476,1.238996021017181) -- (-1.349862258692196,1.0746556467304906);
\draw (-1.3025099474570476,1.238996021017181) -- (-1.450137741307804,1.3253443532695104);
\draw (-0.4,1.6)-- (-1.6,0.);
\draw (-1.063,0.716) -- (-1.108,0.881);
\draw (-1.063,0.716) -- (-0.892,0.719);
\draw (-0.4,1.6)-- (1.2,0.4);
\draw (0.484,0.937) -- (0.319,0.892);
\draw (0.484,0.937) -- (0.481,1.108);
\draw (1.2,2.8)-- (-0.4,1.6);
\draw (0.316,2.137) -- (0.319,2.308);
\draw (0.316,2.137) -- (0.481,2.092);
\draw [->] (5.2,1.6) -- (2.8,1.6);
\draw (9.6,-0.4)-- (9.6,1.6);
\draw (9.6,0.705) -- (9.735,0.6);
\draw (9.6,0.705) -- (9.465,0.6);
\draw (9.6,1.6)-- (11.2,0.4);
\draw (10.484,0.937) -- (10.319,0.892);
\draw (10.484,0.937) -- (10.481,1.108);
\draw (11.2,2.8)-- (9.6,1.6);
\draw (10.316,2.137) -- (10.319,2.308);
\draw (10.316,2.137) -- (10.481,2.092);
\draw (8.8,1.6)-- (8.8,3.6);
\draw (8.8,2.705) -- (8.935,2.6);
\draw (8.8,2.705) -- (8.665,2.6);
\draw (7.6,3.2)-- (8.8,1.6);
\draw (8.263,2.316) -- (8.092,2.319);
\draw (8.263,2.316) -- (8.308,2.481);
\draw (8.8,1.6)-- (6.8,2.4);
\draw (7.702509947457047,2.038996021017181) -- (7.850137741307804,2.12534435326951);
\draw (7.702509947457047,2.038996021017181) -- (7.749862258692194,1.8746556467304898);
\draw (6.8,0.8)-- (8.8,1.6);
\draw (7.897490052542953,1.238996021017181) -- (7.8501377413078055,1.0746556467304906);
\draw (7.897490052542953,1.238996021017181) -- (7.749862258692197,1.3253443532695104);
\draw (8.8,1.6)-- (7.6,0.);
\draw (8.137,0.716) -- (8.092,0.881);
\draw (8.137,0.716) -- (8.308,0.719);
\begin{scriptsize}
\draw [fill=qqqqff] (-0.4,-0.4) circle (2.5pt);
\draw [fill=qqqqff] (-0.4,1.6) circle (2.5pt);
\draw[color=black] (-0.02,1.6) node {$v$};
\draw[color=black] (-0.02,0.78) node {$e_1$};
\draw [fill=qqqqff] (-0.4,3.6) circle (2.5pt);
\draw[color=black] (-0.02,2.78) node {$e_2$};
\draw[color=black] (4.06,1.9) node {$\pi$};
\draw [fill=qqqqff] (9.6,-0.4) circle (2.5pt);
\draw [fill=qqqqff] (9.6,1.6) circle (2.5pt);
\draw[color=black] (9.6,2) node {$\tv_1$};
\draw[color=black] (10.02,0.84) node {$\widetilde{e}_1$};
\draw [fill=qqqqff] (8.8,1.6) circle (2.5pt);
\draw[color=black] (8.9,1.2) node {$\tv_2$}; 
\draw [fill=qqqqff] (8.8,3.6) circle (2.5pt);
\draw[color=black] (9.22,2.84) node {$\widetilde{e}_2$};
\end{scriptsize}
\end{tikzpicture}
\end{center}
\caption{}
\label{Fig: 3-L}
\end{figure}

\noindent

However this is contradiction to the fact that $\tGamma$ is hyperbolic, because $\tv_{1},\tv_{2}$ have odd valency.
As a consequence of the claim we may assume $\tC_{1}=C_{1}$ and $\tC_{2}$ has no edge.

We check the second step $\tC_{2}=\emptyset$.
Suppose it is nonempty, i.e., $\tC_{2}$ consists of only vertices.
By genus induction $\tGamma$ is a cactus boundary so that it is depicted as Fig.\ \ref{Fig: 3-M}.
\begin{figure}[ht]
\begin{center}
\begin{tikzpicture}[line cap=round,line join=round,>=triangle 45,x=.7cm,y=.7cm]
\draw (-0.4,1.6)-- (0.4,2.8);
\draw (0.05824352060364902,2.2873652809054743) -- (0.11232678973560895,2.1251154735095943);
\draw (0.05824352060364902,2.2873652809054743) -- (-0.11232678973560895,2.2748845264904065);
\draw (0.4,2.8)-- (1.6,2.8);
\draw (1.105,2.8) -- (1.,2.665);
\draw (1.105,2.8) -- (1.,2.935);
\draw (1.6,2.8)-- (2.4,1.6);
\draw (2.0582435206036496,2.1126347190945265) -- (1.8876732102643916,2.1251154735095943);
\draw (2.0582435206036496,2.1126347190945265) -- (2.1123267897356097,2.2748845264904065);
\draw (2.4,1.6)-- (1.6,0.4);
\draw (1.9417564793963504,0.9126347190945263) -- (1.8876732102643905,1.0748845264904061);
\draw (1.9417564793963504,0.9126347190945263) -- (2.1123267897356084,0.9251154735095942);
\draw (1.6,0.4)-- (0.4,0.4);
\draw (0.895,0.4) -- (1.,0.535);
\draw (0.895,0.4) -- (1.,0.265);
\draw (0.4,0.4)-- (-0.4,1.6);
\draw (-0.05824352060364902,1.0873652809054741) -- (0.11232678973560895,1.0748845264904061);
\draw (-0.05824352060364902,1.0873652809054741) -- (-0.11232678973560895,0.9251154735095942);
\draw (2.4,1.6)-- (2.8,2.4);
\draw (2.6469574275274956,2.0939148550549915) -- (2.720747670784989,1.9396261646075064);
\draw (2.6469574275274956,2.0939148550549915) -- (2.479252329215011,2.060373835392495);
\draw (2.8,2.4)-- (3.6,2.4);
\draw (3.305,2.4) -- (3.2,2.265);
\draw (3.305,2.4) -- (3.2,2.535);
\draw (3.6,2.4)-- (4.,1.6);
\draw (3.8469574275274954,1.9060851449450091) -- (3.679252329215011,1.939626164607506);
\draw (3.8469574275274954,1.9060851449450091) -- (3.920747670784989,2.0603738353924945);
\draw (4.,1.6)-- (3.6,0.8);
\draw (3.7530425724725047,1.1060851449450095) -- (3.679252329215011,1.2603738353924947);
\draw (3.7530425724725047,1.1060851449450095) -- (3.920747670784989,1.1396261646075063);
\draw (3.6,0.8)-- (2.8,0.8);
\draw (3.095,0.8) -- (3.2,0.935);
\draw (3.095,0.8) -- (3.2,0.665);
\draw (2.8,0.8)-- (2.4,1.6);
\draw (2.5530425724725045,1.2939148550549915) -- (2.720747670784989,1.2603738353924947);
\draw (2.5530425724725045,1.2939148550549915) -- (2.479252329215011,1.1396261646075063);
\draw (4.,1.6)-- (4.8,2.8);
\draw (4.8,2.8)-- (6.,2.8);
\draw (6.,2.8)-- (6.8,1.6);
\draw (6.8,1.6)-- (6.,0.4);
\draw (6.,0.4)-- (4.8,0.4);
\draw (4.8,0.4)-- (4.,1.6);
\draw [rotate around={0.:(3.2,1.6)}] (3.2,1.6) ellipse (5.122499389946278 and 3.2);
\begin{scriptsize}
\draw [fill=qqqqff] (-0.4,1.6) circle (2.5pt);
\draw [fill=qqqqff] (0.4,2.8) circle (2.5pt);
\draw [fill=qqqqff] (1.6,2.8) circle (2.5pt);
\draw[color=black] (1.06,2.06) node {$C_1$};
\draw [fill=qqqqff] (2.4,1.6) circle (2.5pt);
\draw [fill=qqqqff] (1.6,0.4) circle (2.5pt);
\draw [fill=qqqqff] (0.4,0.4) circle (2.5pt);
\draw [fill=qqqqff] (2.8,2.4) circle (2.5pt);
\draw [fill=qqqqff] (3.6,2.4) circle (2.5pt);
\draw [fill=qqqqff] (4.,1.6) circle (2.5pt);
\draw [fill=qqqqff] (3.6,0.8) circle (2.5pt);
\draw [fill=qqqqff] (2.8,0.8) circle (2.5pt);
\draw[color=black] (5.46,2.06) node {$C_2$};
\end{scriptsize}
\end{tikzpicture}
\end{center}
\caption{}
\label{Fig: 3-M}
\end{figure}

However in Fig.\ \ref{Fig: 3-M} it is impossible to construct $\tGamma_{\sr}$ because $\tGamma_\sr$, if any, should intersect the puncture $C_{2}$ as in Fig.\ \ref{Fig: 3-m}.
\begin{figure}[ht]
\begin{center}
\begin{tikzpicture}[line cap=round,line join=round,>=triangle 45,x=.7cm,y=.7cm]
\draw[color=gray,arrow data={.25}{stealth},arrow data={.5}{stealth},arrow data={.75}{stealth},arrow data={.99}{stealth}] plot [smooth, tension=.3pt] coordinates
{(1.,0.2) (0.2,0.2)   (-0.6,1.4)   (-0.615333333333333,1.8273333333333346)   (0.2013333333333335,3.0273333333333343)   (1.8,3.)   (2.4,2.2)   (2.8,2.8)   (3.768,2.8106666666666675)   (4.4,1.6)   (4.084666666666666,1.0106666666666682) (3.768,0.41066666666666823)   (2.7513333333333327,0.41066666666666823)   (2.3846666666666665,1.0606666666666682)   (1.768,0.21066666666666828) (1.,0.2)};
\draw (-0.4,1.6)-- (0.4,2.8);
\draw (0.048536267169707506,2.272804400754561) -- (0.09360565811300743,2.137596227924661);
\draw (0.048536267169707506,2.272804400754561) -- (-0.09360565811300743,2.262403772075338);
\draw (0.4,2.8)-- (1.6,2.8);
\draw (1.0875,2.8) -- (1.,2.6875);
\draw (1.0875,2.8) -- (1.,2.9125);
\draw (1.6,2.8)-- (2.4,1.6);
\draw (2.048536267169708,2.127195599245438) -- (1.906394341886993,2.137596227924661);
\draw (2.048536267169708,2.127195599245438) -- (2.093605658113008,2.262403772075338);
\draw (2.4,1.6)-- (1.6,0.4);
\draw (1.951463732830292,0.9271955992454383) -- (1.9063943418869922,1.062403772075338);
\draw (1.951463732830292,0.9271955992454383) -- (2.093605658113007,0.9375962279246615);
\draw (1.6,0.4)-- (0.4,0.4);
\draw (0.9125,0.4) -- (1.,0.5125);
\draw (0.9125,0.4) -- (1.,0.2875);
\draw (0.4,0.4)-- (-0.4,1.6);
\draw (-0.048536267169707506,1.0728044007545612) -- (0.09360565811300743,1.062403772075338);
\draw (-0.048536267169707506,1.0728044007545612) -- (-0.09360565811300743,0.9375962279246615);
\draw (2.4,1.6)-- (2.8,2.4);
\draw (2.639131189606246,2.0782623792124917) -- (2.7006230589874907,1.9496884705062536);
\draw (2.639131189606246,2.0782623792124917) -- (2.499376941012509,2.0503115294937446);
\draw (2.8,2.4)-- (3.6,2.4);
\draw (3.2875,2.4) -- (3.2,2.2875);
\draw (3.2875,2.4) -- (3.2,2.5125);
\draw (3.6,2.4)-- (4.,1.6);
\draw (3.8391311896062468,1.9217376207875072) -- (3.6993769410125097,1.9496884705062547);
\draw (3.8391311896062468,1.9217376207875072) -- (3.9006230589874913,2.0503115294937455);
\draw (4.,1.6)-- (3.6,0.8);
\draw (3.7608688103937546,1.1217376207875074) -- (3.6993769410125097,1.2503115294937457);
\draw (3.7608688103937546,1.1217376207875074) -- (3.9006230589874913,1.149688470506255);
\draw (3.6,0.8)-- (2.8,0.8);
\draw (3.1125,0.8) -- (3.2,0.9125);
\draw (3.1125,0.8) -- (3.2,0.6875);
\draw (2.8,0.8)-- (2.4,1.6);
\draw (2.5608688103937536,1.278262379212492) -- (2.7006230589874907,1.2503115294937446);
\draw (2.5608688103937536,1.278262379212492) -- (2.499376941012509,1.1496884705062538);
\draw (4.,1.6)-- (4.8,2.8);
\draw (4.8,2.8)-- (6.,2.8);
\draw (6.,2.8)-- (6.8,1.6);
\draw (6.8,1.6)-- (6.,0.4);
\draw (6.,0.4)-- (4.8,0.4);
\draw (4.8,0.4)-- (4.,1.6);
\draw [rotate around={0.:(3.2,1.6)}] (3.2,1.6) ellipse (5.122499389946278 and 3.2);
\begin{scriptsize}
\draw [fill=qqqqff] (-0.4,1.6) circle (2.5pt);
\draw [fill=qqqqff] (0.4,2.8) circle (2.5pt);
\draw [fill=qqqqff] (1.6,2.8) circle (2.5pt);
\draw [fill=qqqqff] (2.4,1.6) circle (2.5pt);
\draw [fill=qqqqff] (1.6,0.4) circle (2.5pt);
\draw [fill=qqqqff] (0.4,0.4) circle (2.5pt);
\draw [fill=qqqqff] (2.8,2.4) circle (2.5pt);
\draw [fill=qqqqff] (3.6,2.4) circle (2.5pt);
\draw [fill=qqqqff] (4.,1.6) circle (2.5pt);
\draw [fill=qqqqff] (3.6,0.8) circle (2.5pt);
\draw [fill=qqqqff] (2.8,0.8) circle (2.5pt);
\end{scriptsize}
\end{tikzpicture}
\end{center}
\caption{}
\label{Fig: 3-m}
\end{figure}

Now we get $\tC_{2}=\emptyset$ (and $\tC_{1}=C_{1}$).
Hence we obtain a semi-ribbon graph $(\tS,\tGamma,\tGamma_{\sr})$ such that $\tGamma_{\sr}$ lies outside the puncture $C_{2}$.
Therefore $(\tS,\tGamma,\tGamma_{\sr})$ satisfies the theorem by genus induction.
Since $S$ is a connected sum of $\tS$ and a cylinder along $C_{1},C_{2}$, $(S,\Gamma,\Gamma_{\sr})$ also satisfies the theorem.

The statement on the number of full Eulerian paths from semi-ribbon graphs is clear because $\Gamma$ is a cactus boundary.
This completes the proof of the theorem.


\section{$\omega$-limit sets}\label{sec: omega}

In this section we prove Theorem \ref{th: main: analytic version}.
This will be done by generalizing the theorem to Theorem \ref{th: main: omega} using $\cC^{1}$-flows on an oriented $\cC^{1}$-surface with corner.
However we impose further local condition on the $\cC^{1}$-flows, namely analytic sector decomposition (Definition \ref{def: analytic sector decomposition}).

Let $(S,\Phi)$ be a $\cC^1$-flow on an oriented $\cC^{1}$-surface $S$ with corner.
In this section we always set $\Gamma$ to be the corresponding graph to an $\omega$-limit set $\omega(x)$ if $\omega(x)$ is a polycycle.
Recall that this correspondence is given by 
	$$
	\begin{aligned}
	&
	\mbox{free orbits in $\omega(x)$ $\leftrightarrow$ edges of $\Gamma$,} 
	\\
	&
	\mbox{singular points in $\omega(x)$ $\leftrightarrow$ vertices of $\Gamma$}.
	\end{aligned}
	$$

In \S\ref{subsec: genus 0} we explain the classification of $(S,\Gamma)$ in \cite{JL} when $S$ has genus $0$.
This result will serve as the initial genus induction hypothesis for the higher genus.
In \S\ref{subsec: higher genus} we state Theorem \ref{th: main: omega} and prepare some lemmas towards the proof of the theorem. 
This theorem is a generalization of Theorem \ref{th: main: analytic version}.
In \S\ref{subsec: Proof of Theorem omega} we prove Theorem \ref{th: main: omega} using genus induction.


\subsection{Genus 0}\label{subsec: genus 0}

It is known that there exists an open neighborhood $U$ of an isolated singular point $p$ of an analytic surface flow $\Phi$ on an analytic surface satisfying \textit{either}
	\begin{itemize}
	\item[(i)]
	the orbits of $\Phi|_{U}$ are concentric periodic orbits \textit{or}
	\item[(ii)]
	$U$ is decomposed into finitely many sectors centered at $p$ each of which consists of either hyperbolic, repelling, contracting or elliptic orbits (ref.\ \cite[pp.85--86]{AI}).
	\end{itemize}

In the case when an analytic surface $S$ has corner, we have a similar decomposition. 
If an isolated singularity $p$ lies in the boundary of $S$, the adjacent boundary sectors are required to be separatices of sector decomposition.

\begin{defn}\label{def: analytic sector decomposition}
A $\cC^{1}$-flow $\Phi$ has \textit{analytic sector decomposition} if at any singular point, there is a local  $\cC^{1}$-diffeomorphism with an analytic flow.
\end{defn}

This property amounts to that $\Phi$ has only isolated singularity and allows a sector decomposition at a singular point as (i) and (ii) up to $\cC^{1}$-diffeomorphism.
It is automatic that if $\Phi$ has analytic sector decomposition, the set of singularities $\Sing(\Phi)$ is finite because $S$ is compact and thus the sector decomposition cannot be accumulated.

For an oriented surface $S$ with corner, we denote by $S^{\op}$ the orientation-reversed surface.
The $\omega$-limit sets for the genus $0$ are described in terms of semi-ribbon graphs:

\begin{thm}\label{th: JL}
Suppose that $g(S)=0$ and that $\Phi$ has  analytic sector decomposition.
For a given $x\in S$, if $x$ is not a singular point nor a periodic point, $\Gamma$ is a cactus boundary.

Moreover  there is a semi-ribbon graph $\Gamma_{\sr}$ of $\Gamma$ in either $S$ or $S^{\op}$ such that $\Phi(t,x)\in \Gamma_{\sr}\setminus\Gamma$ for all sufficiently large $t$ and  $\lim_{t\to\infty}d(\Phi(t,y),\omega(x))=0$ for any $y\in \Gamma_{\sr}$.
\end{thm}

\proof
Suppose first $S$ has no boundary, i.e., $S=S^{2}$.
Then the theorem is only restatement of \cite[Theorem C]{JL}.
In fact \textit{loc}.\ \textit{cit}.\ is stated for analytic flows on $S^{2}$.
But the proof uses only analytic sector decomposition.
We omit the full details, but refer to the remark after the proof.

We will define $\Gamma_{\sr}$.
We follow the idea in the proof of \cite[Lemma 4.8]{JL}.
Let $x_{0}$ be any smooth point on $\omega(x)$ and $I$ be an embedded closed interval in $S$ whose one end is $x_{0}$ and the other one lies outside $\omega(x)$.
Then $\Phi(t,x)$ intersects $I$ for some $t=t_{0}$.
The Poincar\'e return map sends $\Phi(t_{0},x)$ to $\Phi(t_{1},x)$ for some $t_{1}>t_{0}$.
Now we define $\Gamma_{\sr}$ as the strip bounded by $\Phi(t,x)$ $(t_{0}\le t\le t_{1})$ and the subinterval in $I$ connecting $\Phi(t_{0},x)$ and $\Phi(t_{1},x)$.
This $\Gamma_{\sr}$ is a semi-ribbon graph of $(S,\Gamma)$ if it has the clockwise orientation.
Otherwise $(S^{\op},\Gamma,\Gamma_{\sr})$ is a semi-ribbon graph.

In the case when $S$ has boundary, it is not difficult to extend $(S,\Phi)$ to $(S_{\cl},\Phi_{\cl})$ such that $\Phi_{\cl}$ has analytic sector decomposition.
This extension appears after this paragraph. 
Now since $S_{\cl}=S^{2}$, the above argument shows the theorem.

We fix $C$ any boundary circle of $S$ and $D$ be the capping-off disk.
We need to extend the flow $\Phi$ on $C$ to the whole disk $D$ with analytic sector decomposition. 
We may assume that $D$ is the unit disk in $\rr^{2}$.
If $C$ has only one-way direction flow (possibly with singularities), one can extend this to $D$ with the concentric orbits centered at the origin $0$.
It is easy to see that this extension flow satisfies the analytic sector decomposition near $D$.  

\begin{figure}[ht]
\begin{center}
\begin{tikzpicture}[line cap=round,line join=round,>=triangle 45,x=0.7cm,y=0.7cm]

\draw [arrow data={.5}{stealth},shift={(1.,0.4)}] plot[domain=-1.5707963267948966:0.,variable=\t]({1.*3.2*cos(\t r)+0.*3.2*sin(\t r)},{0.*3.2*cos(\t r)+1.*3.2*sin(\t r)});
\draw [arrow data={.5}{stealth},shift={(1.,0.4)}] plot[domain=0.:1.5707963267948966,variable=\t]({1.*3.2*cos(\t r)+0.*3.2*sin(\t r)},{0.*3.2*cos(\t r)+1.*3.2*sin(\t r)});
\draw [arrow data={.5}{stealth},shift={(1.,0.4)}] plot[domain=2.351789232345534:1.5707963267948966,variable=\t]({1.*3.2*cos(\t r)+0.*3.2*sin(\t r)},{0.*3.2*cos(\t r)+1.*3.2*sin(\t r)});
\draw [arrow data={.5}{stealth},shift={(1.,0.4)}] plot[domain=2.351789232345534:3.141592653589793,variable=\t]({1.*3.2102959365142647*cos(\t r)+0.*3.2102959365142647*sin(\t r)},{0.*3.2102959365142647*cos(\t r)+1.*3.2102959365142647*sin(\t r)});
\draw [arrow data={.5}{stealth},shift={(1.,0.4)}] plot[domain=4.71238898038469:3.141592653589793,variable=\t]({1.*3.2*cos(\t r)+0.*3.2*sin(\t r)},{0.*3.2*cos(\t r)+1.*3.2*sin(\t r)});
\draw [arrow data={.5}{stealth},shift={(10.2,0.4)}] plot[domain=-1.5707963267948966:0.,variable=\t]({1.*3.2*cos(\t r)+0.*3.2*sin(\t r)},{0.*3.2*cos(\t r)+1.*3.2*sin(\t r)});
\draw [arrow data={.5}{stealth},shift={(10.2,0.4)}] plot[domain=0.:1.5707963267948966,variable=\t]({1.*3.2*cos(\t r)+0.*3.2*sin(\t r)},{0.*3.2*cos(\t r)+1.*3.2*sin(\t r)});
\draw [arrow data={.5}{stealth},shift={(10.2,0.4)}] plot[domain=2.3517892323455336:1.5707963267948966,variable=\t]({1.*3.2*cos(\t r)+0.*3.2*sin(\t r)},{0.*3.2*cos(\t r)+1.*3.2*sin(\t r)});
\draw [arrow data={.5}{stealth},shift={(10.2,0.4)}] plot[domain=2.3517892323455336:3.141592653589793,variable=\t]({1.*3.2102959365142643*cos(\t r)+0.*3.2102959365142643*sin(\t r)},{0.*3.2102959365142643*cos(\t r)+1.*3.2102959365142643*sin(\t r)});
\draw [arrow data={.5}{stealth},shift={(10.2,0.4)}] plot[domain=4.71238898038469:3.141592653589793,variable=\t]({1.*3.2*cos(\t r)+0.*3.2*sin(\t r)},{0.*3.2*cos(\t r)+1.*3.2*sin(\t r)});
\draw[arrow data={.5}{stealth}] (10.2,3.6)-- (10.2,0.4);
\draw[arrow data={.5}{stealth}] (10.2,0.4)-- (10.2,-2.8);
\draw[arrow data={.5}{stealth}] (10.2,0.4)-- (7.94,2.68);
\draw[arrow data={.5}{stealth}] (7.,0.4)-- (10.2,0.4);
\begin{scriptsize}
\draw [fill=qqqqff] (1.,-2.8) circle (2.5pt);
\draw [fill=qqqqff] (4.2,0.4) circle (2.5pt);
\draw [fill=qqqqff] (1.,3.6) circle (2.5pt);
\draw [fill=qqqqff] (-1.26,2.68) circle (2.5pt);
\draw [fill=qqqqff] (-2.2,0.4) circle (2.5pt);
\draw [fill=qqqqff] (10.2,0.4) circle (2.5pt);
\draw [fill=qqqqff] (10.2,-2.8) circle (2.5pt);
\draw [fill=qqqqff] (13.4,0.4) circle (2.5pt);
\draw [fill=qqqqff] (10.2,3.6) circle (2.5pt);
\draw [fill=qqqqff] (7.94,2.68) circle (2.5pt);
\draw [fill=qqqqff] (7.,0.4) circle (2.5pt);
\end{scriptsize}
\end{tikzpicture}
\end{center}
\caption{extension of the flow on $C$ to the boundaries of sectors}
\label{Fig: 4-alpha}
\end{figure}
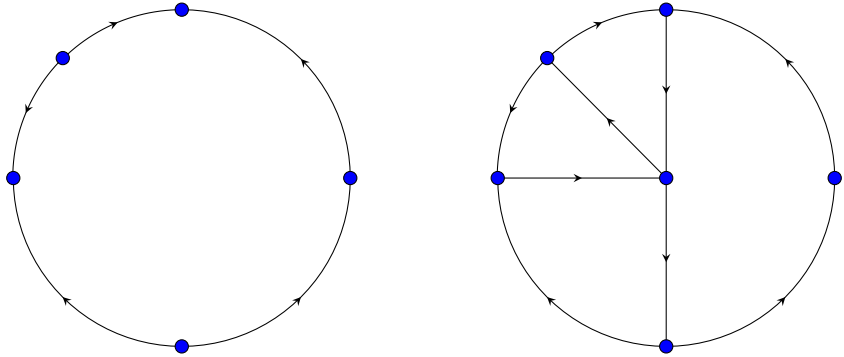

Otherwise there is a decomposition of $C$ by the non-overlapping proper arcs such that any two adjacent arcs have opposite flow directions.
Thus there is the corresponding sector decomposition of $D$ in a way that each sector has the above arc.  
We extend the flow $\Phi|_{C}$ to the union of the boundaries of the sectors such that each boundary has one-way flow direction.
See Fig.\ \ref{Fig: 4-alpha}.

Now we extend the flow on the boundary of each sector to the whole sector respectively as follows.
Let $T$ be any sector among the above. 
Let $(T,\Bd(T))\approx (D,C)$ be any homeomorphism such that the interiors of $T,D$ are diffeomorphic and $\Bd(T),C$ are piecewise diffeomorphic.
We endow a flow $C$ by the flow on $\Bd(T)$ using this homeomorphism.
We can extend this flow on $C$ to the whole $D$ concentrically as before.
Now the induced extension flow on $T$ satisfies the analytic sector decomposition near $T$.
Hence so does the extension flow on $D$.

By applying the above procedure to each capping-off disk, we obtain the extension flow with analytic sector decomposition on $S_{\cl}=S^{2}$.
\qed\vskip.3cm

Note that Lemma \ref{lem: geneus 0}  is the graph-theoretic counterpart of this theorem.

\begin{rk}\label{rk: gap}
A gap in the proof of \cite[Theorem C]{JL} is fixed in \cite{BL} while the statement itself remains true.
This is due to an erratum in the second assertion of \cite[Lemma 4.6]{JL}, namely ``equal orientation'' of semi-open flow boxes.
To be precise the proof of \cite[Lemma 4.6]{JL} is divided into the two sub-cases, isolated singularities and non-isolated ones; an erratum occurs in the latter case.
Under our assumption i.e., the finite singularity assumption, the assertion of \cite[Lemma 4.6]{JL} holds without change, so that our argument in the above proof regarding \cite[Theorem C]{JL} remains valid.
\end{rk}


\subsection{Higher genus}\label{subsec: higher genus}

We give the topological classification theorem of $(S,\Phi)$, where $S$ is an oriented $\cC^{1}$-surface with corner and $\Phi$ is a $\cC^{1}$-flow with analytic sector decomposition.

\begin{thm}\label{th: main: omega}
Suppose that $g(S)\ge0$ and that $\Phi$ has analytic sector decomposition and that $\omega(x)$ is a polycycle.
Then there is a semi-ribbon graph $\Gamma_{\sr}$ of $\Gamma$ in either $S$ or $S^{\op}$ such that 
	\begin{itemize}
	\item[(i)]
	$\Phi(t,x)\in \Gamma_{\sr}\setminus\Gamma$ for all sufficiently large $t$,
	\item[(ii)]  
	$\lim_{t\to \infty}d(\Phi(t,y),\omega(x))=0$ for any  $y\in \Gamma_{\sr}$.
	\end{itemize}
Hence either $(S,\Gamma,\Gamma_{\sr})$ or $(S^{\op},\Gamma,\Gamma_{\sr})$ is described as in Theorem \ref{th: number of Eulerian paths}.

Conversely for any semi-ribbon graph $(S,\Gamma,\Gamma_{\sr})$, there is a flow $\Phi$ with analytic sector decomposition and $x\in S$ such that $\omega(x)=\Gamma$ and $\Phi(t,x)$ satisfying (i) and (ii) in the above.
\end{thm}

The proof of the theorem appears in \S\ref{subsec: Proof of Theorem omega}.
This theorem generalizes Theorem \ref{th: main: analytic version} which classifies the analytic flows on oriented analytic surfaces.
And it has a graph-theoretic counterpart of Theorem \ref{th: number of Eulerian paths}.

\begin{rk}\label{rk: polycycle}
Note that the polycycle assumption on $\omega(x)$ in Theorem \ref{th: main: omega} does not appear in Theorem \ref{th: JL}. 
For, in the genus $0$ case, $\omega(x)$ always turns out to be a polycycle, unless $x$ is a singular point or a periodic point.
\end{rk}

Now we need preliminary steps towards the proof of Theorem \ref{th: main: omega}.
These steps are not for the converse direction of the theorem, as it is rather easy to prove.

Note that any edge $e$ of $\Gamma$ is a separatrix of a sector.
Let $v$ be either $t(e)$ or $h(e)$ (i.e., $e$ is an \textit{adjoining} edge of $v$).
Then there are precisely two adjacent sectors near $v$ having $e$ as a common separatrix.

\begin{lem}\label{lem: condition}
$(1)$
Every edge $e$ of $\Gamma$ is a separatrix of at least one hyperbolic sector centered at $v$.

$(2)$
There is a hyperbolic sector in the above item $(1)$ containing $\Phi(t_{k},x)$ in its interior for some sequence $(t_{k})$ with $t_{k}\to\infty$ as $k\to\infty$.
\end{lem}

\proof
(1)
We need to exclude the following possibility: each of the two adjacent sectors centered at $v$ is one of repelling, attracting, elliptic sectors. 
There are six pairs of types of adjacent sectors.
Suppose the pair is elliptic-elliptic.
We may assume $x$ is contained in the interior of one elliptic sector since if $x\in e$ then $\omega(x)$ is a one-point set and thus $e$ cannot be a part of $\omega(x)$.
Even under this assumption we have $\omega(x)=\{v\}$.
This is contradiction.
For the other pairs, we meet a similar contradiction to that $e\subset \omega(x)$.
We omit the details.

(2)
There are only two possibilities: (i) $e$ is a separatrix of two adjacent hyperbolic sectors, (ii) $e$ is a separatrix of only one hyperbolic sector.
In the case (i), we need to check $\Phi(t_{k},x)$ lies in the interior of at least one hyperbolic sector for infinitely many $t_{k}\in\rr$ with $t_{k}\to \infty$ as $k\to \infty$.
But this is obvious because otherwise, any smooth point in $e$ cannot be an element of $\omega(x)$.

In the case (ii), $\Phi(t,x)$ cannot lie in the non-hyperbolic sector by the proof of the item (1) for any $t$.
So by the discussion in the proof of (i), $\Phi(t_{k},x)$ lies in the interior of  the hyperbolic sector for infinitely many $t_{k}\in\rr$ with $t_{k}\to \infty$ as $k\to\infty$. 
This completes the proof of the lemma.
\qed\vskip.3cm

The above lemma does not assert $\Gamma$ is hyperbolic.
Thus we do not know if $\Gamma$ admits a full Eulerian path at this moment.
However a similar but much simpler argument used in the proof of Lemma \ref{lem: existence of minimal closed path}, we deduce the following:

\begin{lem}\label{lem: existence of simple closed path}
There exists a simple closed path in $\Gamma$.
\end{lem}

\proof
We take any vertex $v_{0}$ as a starting point and then a simple path, say $\gamma$, from $v_{0}$ with the maximal number of edges.
If the terminal vertex of $\gamma$, say $v_{1}$, becomes $v_{0}$,  then we are done.
Otherwise there is an edge in $\Gamma$ with the tail vertex $v_{1}$ by the definition of polycycle which is not contained in $\gamma$. 
Thus we obtain a longer path, say $\tgamma$ than $\gamma$, but this contradicts the maximality unless $\tgamma$ is not simple.
Now $\tgamma$ contains a simple closed subpath containing $v_{1}$.
\qed\vskip.3cm


\subsection{Proof of Theorem \ref{th: main: omega}.}\label{subsec: Proof of Theorem omega}
The idea goes as in the proof of Theorem \ref{th: number of Eulerian paths}.
We will show that $\Gamma$ admits a semi-ribbon graph $\Gamma_{\sr}$ so that the topological description of $(S,\Gamma)$ will follow from Theorem \ref{th: number of Eulerian paths}.
And then we will check the assertions on convergence in the theorem.

First we know from Lemma \ref{lem: existence of simple closed path} that there is at least one simple closed path of $\Gamma$.

Suppose there is a simple closed path of $\Gamma$ with the trivial homology class in $H^{1}(S,\zz)$.
It is an embedded disk boundary $\Bd(D)$.
By Lemma \ref{lem: existence of minimal closed path} we may assume $\Int(D)$ contains no edge.
Here we notice that the assumption of Lemma \ref{lem: existence of minimal closed path} is satisfied by the definition of polycycle.
If $\Gamma=\Bd(D)$, then $\Gamma$ is contained in the genus $0$ surface.
If $x$ lies in $D$,  Theorem \ref{th: JL} implies Theorem \ref{th: main: omega}.
If $x$ lies outside $D$, we set $\Gamma_{\sr}$ to be the strip outside $D$ as in the proof of Theorem \ref{th: JL}.

We assume $\Gamma\neq\Bd(D)$ so that $\Gamma$ has an edge outside $D$.
We will use the induction on $\#E$.
The initial induction hypothesis is the case when $\#E=1$, i.e., $\Gamma$ is a loop with one vertex.
If the loop encloses an embedded disk in $S$, we are already done.
Otherwise we need to proceed to the next case when there is no nomologically trivial simple closed path in $\Gamma$.  

Let us consider the shrinking of $D$ as in the proof of Lemma \ref{lem: homology 0}.
Let $(\oS,\oPhi,\oGamma)$ be the induced triple of $(S,\Phi,\Gamma)$ from shrinking.
We claim $\oPhi$ has also analytic sector decomposition.
The change occurs only near at the shrunken vertex of $\oGamma$.
We take vertices $v,v'\in\Bd(D)$ and edges $e,e'$ outside $D$ adjoining $v,v'$.
Let $P$ be the path in $\Bd(D)$ connecting $v,v'$.
See Fig.\ \ref{Fig: 4-A} (in this figure, $P$ is the path from $v$ to $v'$).
\begin{figure}[ht]
\begin{center}
\begin{tikzpicture}[line cap=round,line join=round,>=triangle 45,x=1cm,y=1cm]
\draw (0.,2.4)-- (2.4,2.4);
\draw (1.305,2.4) -- (1.2,2.265);
\draw (1.305,2.4) -- (1.2,2.535);
\draw (2.4,2.4)-- (2.4,0.);
\draw (2.4,1.095) -- (2.265,1.2);
\draw (2.4,1.095) -- (2.535,1.2);
\draw (2.4,0.)-- (0.,0.);
\draw (1.095,0.) -- (1.2,0.135);
\draw (1.095,0.) -- (1.2,-0.135);
\draw (0.,0.)-- (0.,2.4);
\draw (0.,1.305) -- (0.135,1.2);
\draw (0.,1.305) -- (-0.135,1.2);
\draw (3.6,2.4)-- (2.4,2.4);
\draw (2.895,2.4) -- (3.,2.535);
\draw (2.895,2.4) -- (3.,2.265);
\draw (2.4,2.4)-- (2.4,3.6);
\draw (2.4,3.105) -- (2.535,3.);
\draw (2.4,3.105) -- (2.265,3.);
\draw (-1.2,0.)-- (0.,0.);
\draw (-0.495,0.) -- (-0.6,-0.135);
\draw (-0.495,0.) -- (-0.6,0.135);
\draw (0.,0.)-- (0.,-1.2);
\draw (0.,-0.705) -- (-0.135,-0.6);
\draw (0.,-0.705) -- (0.135,-0.6);
\draw [->,decorate,decoration=zigzag] (4.8,1.2) -- (7.2,1.2);
\draw (8.4,1.2)-- (9.6,1.2);
\draw (9.105,1.2) -- (9.,1.065);
\draw (9.105,1.2) -- (9.,1.335);
\draw (9.6,1.2)-- (9.6,0.);
\draw (9.6,0.495) -- (9.465,0.6);
\draw (9.6,0.495) -- (9.735,0.6);
\draw (10.8,1.2)-- (9.6,1.2);
\draw (10.095,1.2) -- (10.2,1.335);
\draw (10.095,1.2) -- (10.2,1.065);
\draw (9.6,1.2)-- (9.6,2.4);
\draw (9.6,1.905) -- (9.735,1.8);
\draw (9.6,1.905) -- (9.465,1.8);
\draw[color=gray,arrow data={.25}{stealth},arrow data={.5}{stealth},arrow data={.75}{stealth}] plot [smooth, tension=.4pt] coordinates
{ (3.6,2.)  (3.,1.8) (2.8,0.)  (2.4,-0.4) (0.64,-0.62)  (0.4,-1.2)};
\draw[color=gray,arrow data={.25}{stealth},arrow data={.5}{stealth},arrow data={.75}{stealth}] plot [smooth, tension=.4pt] coordinates
{(-1.2,0.4)  (-0.6,0.6) (-0.4,2.4) (0.,2.8) (1.8,3.) (2.,3.6)};
\draw[color=gray,arrow data={.65}{stealth}] plot [smooth, tension=.4pt] coordinates
{(10.8,0.8) (10.2,0.6) (10.,0.)};
\draw[color=gray,arrow data={.65}{stealth}] plot [smooth, tension=.4pt] coordinates
{(8.4,1.6) (9.04,1.82) (9.2,2.4)};
\begin{scriptsize}
\draw [fill=qqqqff] (0.,2.4) circle (2.5pt);
\draw [fill=qqqqff] (2.4,2.4) circle (2.5pt);
\draw [fill=qqqqff] (2.4,0.) circle (2.5pt);
\draw [fill=qqqqff] (0.,0.) circle (2.5pt);
\draw [fill=qqqqff] (9.6,1.2) circle (2.5pt);
\draw [fill=qqqqff] (1.2,1.2) node {$D$};
\draw [fill=qqqqff] (2.7,2.7) node {$v$};
\draw [fill=qqqqff] (3.8,2.7) node {$e$};
\draw [fill=qqqqff] (-0.3,-0.3) node {$v'$};
\draw [fill=qqqqff] (-0.3,-1.3) node {$e'$};
\draw [color=black] (3.6,1.4) node {$R_1$};
\draw [color=black] (1.2,-1.2) node {$R_2$};
\draw [color=black] (10.8,0.) node {$R$};
\end{scriptsize}
\end{tikzpicture}
\end{center}
\caption{}
\label{Fig: 4-A}
\end{figure}

We may further assume that 
	\begin{enumerate}
	\item
	any other vertices in $P$ than $v,v'$ have no adjoining edges outside $D$,
	\item
	there is a sector centered at $v$ whose separatices are $e,P$,
	\item
	there is a sector centered at $v'$ whose separatices are $e',P$.
	\end{enumerate}
We denote by $R_{1}$ the sector centered at $v$ with the separatices $e,P$.
See Fig.\ \ref{Fig: 4-A}.
Similarly $R_{2}$ denotes the sector centered at $v'$ with the separatrices $e',P$.
Note that $x$ cannot lie in $D$, since otherwise $\omega(x)$ lies in $D$.
Therefore by applying Lemma \ref{lem: condition} (2) to the separatrix $P$, both $R_{1},R_{2}$ are hyperbolic sectors.
Via the shrinking, $R_{1},R_{2}$ are modified into one hyperbolic sector, say $R$.
This proves the claim.

By the induction on $\#E$, the shrunken triple $(\oS,\oPhi,\oGamma)$ satisfies the theorem. 
Thus there is a semi-ribbon graph $\oGamma_{\sr}$.
From Fig.\ \ref{Fig: 4-A} we see $\oGamma_{\sr}$ induces a semi-ribbon graph $\Gamma_{\sr}$ (the gray curves in the figure stand for the smooth boundaries).
Hence by the proof of Theorem \ref{th: number of Eulerian paths}, we obtain a semi-ribbon graph $(S,\Gamma,\Gamma_{\sr})$.
Note that $x$ lies outside $D$ because the edges outside $D$ are parts of $\omega(x)$.
Thus we may assume $x\in \Gamma_{\sr}\setminus \Gamma$.
Since there is a diffeomorphism $\oGamma_{\sr}\setminus \oGamma\approx \Gamma_{\sr}\setminus \Gamma$, we regard $x\in \oGamma_{\sr}\setminus \oGamma$ and the induction hypothesis assures the convergence assertion in the theorem for $\oS$.
Therefore the convergence assertion for $S$ comes from the diffeomorphism $\oGamma_{\sr}\setminus \oGamma\approx \Gamma_{\sr}\setminus \Gamma$. 

Now we assume that there are only homologically nontrivial simple closed paths of $\Gamma$.
We will use the induction on $g(S)$.
The initial induction hypothesis is the case $g(S)=0$, which was done in Theorem \ref{th: JL}.
We consider cutting surgery $\pi\colon \tS\to S$ along any simple closed path $C$.
Let $\tx:=\pi\inv(x)$.
Here we used $x\notin C$ since otherwise $\omega(x)$ is a one-point set or a periodic point.
Since $\pi$ is a finite map (i.e., $\#\pi\inv(z)<\infty$ for any $z\in S$), $\omega(\tx)$ is also a polycycle.
Let $\tGamma$ be the corresponding graph to $\omega(\tx)$.
By genus induction, $\tGamma$ admits a semi-ribbon graph $\tGamma_{\sr}$.
On the other hand we know from the proof of Theorem \ref{th: number of Eulerian paths}, $\pi(\tGamma_{\sr})$ becomes a semi-ribbon graph $\Gamma_{\sr}$ of $\Gamma$.
The convergence assertion follows from the diffeomorphism $\pi\colon\tGamma_{\sr}\setminus\tGamma\approx  \Gamma_{\sr}\setminus \Gamma$.
This completes the proof of one direction of the theorem.

We prove the converse direction.
According to the connected sum decomposition of $S$ in Theorem \ref{th: number of Eulerian paths}, it suffices to take a flow $\Phi$ on the sphere $S^{2}$ whose $\omega$-limit of a point in $S$ has $\Gamma$ as an $\omega$-limit set.
This can be done even real analytically by the explicit construction in \cite[\S6]{JL}.
Now we extend $\Phi$ to the whole $S$ in $\cC^{1}$.


\appendix


\section{Proof of Proposition \ref{prop: pre semi ribbon}}\label{subsec: if part}

The `only if' part is clear from the definition as was observed in Remark \ref{rk: pre semi ribbon hyp}
 
To prove the `if' part, we assume that $\Gamma$ is hyperbolic.
Our proof is similar to the argument used in the proof of Theorem \ref{th: number of Eulerian paths}, but it is rather simpler because we use $\pi\inv(\Gamma)$ instead of $\tGamma$.

Since there is a full Eulerian path, there is at least one simple closed path of $\Gamma$.
If any simple closed path in $\Gamma$ has class $0$ in $H_{1}(S,\zz)$, by Lemma \ref{lem: existence of minimal closed path} there is an embedded disk $D$ such that $\Bd(D)$ is a simple closed path of $\Gamma$ containing no vertex inside $D$.
By shrinking $D$ as in Fig.\ \ref{Fig: 3-I}, we obtain $(\oS,\oGamma)$.
We observe that $\oGamma$ is a hyperbolic graph. 
By the induction on $\#E$, $\oGamma$ admits a pre-semi-graph $\oGamma_{\psr}$.
We construct a pre-semi-graph $\Gamma_{\psr}$ from $\oGamma_{\psr}$ as in Fig.\ \ref{Fig: 4-A}.

Suppose that there is no homologically trivial simple closed path in $\Gamma$.
Let $C$ be a simple closed path of $\Gamma$.
Let $\pi\colon\tS\to S$ be the cutting surgery along $C$.
By Proposition \ref{prop: }, $\pi\inv(\Gamma)$ is a hyperbolic graph.
Suppose first that $\tS$ is connected. 
By the genus induction, $(\tS,\pi\inv(\Gamma))$ admits a pre-semi-ribbon graph $(\pi\inv(\Gamma))_{\psr}$.
It is not difficult to see that the $\pi$-image of $(\pi\inv(\Gamma))_{\psr}$ defines a pre-semi-ribbon graph of $\Gamma$.
In the case $\tS$ is disconnected, the argument is similar.
We omit the details.


\end{document}